\documentclass[12pt,reqno]{amsart}
\newcommand{\mysection}[1]{\section{#1}
      \setcounter{equation}{0}}

\usepackage{color}

\newtheorem{theorem}{Theorem}[section]
\newtheorem{lemma}[theorem]{Lemma}

\theoremstyle{definition}
\newtheorem{assumption}{Assumption}[section]

\theoremstyle{remark}
\newtheorem{remark}{Remark}[section]

\newcommand\bR{\mathbb{R}}

\newcommand\frK{\mathfrak{K}}

\newcommand\cD{\mathcal{D}}
\newcommand\cE{\mathcal{E}}
\newcommand\cF{\mathcal{F}}

\newcommand\cH{\mathcal{H}}

\newcommand\cK{\mathcal{K}}
\newcommand\cL{\mathcal{L}}

\newcommand\cQ{\mathcal{Q}}

\makeatletter
 \newcommand{\sumstar}%
 {\operatornamewithlimits{\sum@\kern-.2em\raise1ex\hbox{*}}}
 \makeatother

\begin{document}

\title[Derivative estimates]
{Higher order derivative estimates for finite-difference 
schemes}

\author[I. Gy\"ongy]{Istv\'an Gy\"ongy}
\address{School of Mathematics,
University of Edinburgh,
King's  Buildings,
Edinburgh, EH9 3JZ, United Kingdom}
\email{gyongy@maths.ed.ac.uk}

\author[N.  Krylov]{Nicolai Krylov}%
\thanks{The work of the second author was partially supported
by NSF grant DMS-0653121}
\address{127 Vincent Hall, University of Minnesota,
Minneapolis,
       MN, 55455, USA}
\email{krylov@math.umn.edu}


\begin{abstract}
We give sufficient conditions under which 
solutions of 
finite-difference 
 schemes  in the space variable for 
 second order possibly 
degenerate parabolic and elliptic equations admit estimates of 
spatial derivatives up to any given order 
independent of the mesh size.
\end{abstract}

\maketitle

\mysection{Introduction}                   \label{section02.04.06}

This is the second part of a series of papers devoted to 
studying the smoothness of solutions to 
finite difference schemes for parabolic and  
elliptic partial differential equations given on the whole 
$\bR^d$. These equations can 
 degenerate, for example 
be just first order PDEs. 
As in \cite{GK}, the first part of this series, 
we consider a grid in $\bR^d$ and a large class of 
monotone finite difference schemes in the space variable 
$x$ in $\bR^d$.

For each small parameter $h>0$  
the 
given grid is dilated by $h$ and for each 
$x\in\bR^d$ it is shifted 
so that $x$ becomes a mesh point.  
We are interested in 
the smoothness in $x$ of the solution $
u_h(t,x)$ 
of the difference scheme.  
 In \cite{GK} estimates,  independent of $h$, for the first
order derivatives of $u_h$ in $x$  
were obtained under general conditions introduced 
there.
In the present paper we investigate the higher 
order derivatives of $u_h$ in $x$. 
The main results  
give estimates, independent of $h$,  
for the derivatives of $u_h$ in $x$ up to 
any given order $m$. The conditions 
extend those from \cite{GK}. 
Using these results in the
continuation  of this paper we estimate the derivatives of $u_h$ 
in $h$, and that allows us to develop a new method of obtaining the 
power series of $u_h$ in $h$. Hence we get accelerated finite-difference 
schemes by using Richardson's extrapolation. 
Namely, under general conditions we show that the 
accuracy of  finite difference schemes for parabolic and 
 elliptic PDEs can 
be improved to any order by taking suitable linear combinations of 
finite difference approximations with different mesh-sizes. 
For elliptic PDEs this result is announced by Theorem 2.4 in \cite{GK} 
and for parabolic PDEs by Theorem \ref{theorem 1.08.02.08} below. 
We hope to develop these results in domains for uniformly
nondegenerate equations later. 

Derivative estimates for finite-difference approximations 
for linear and for nonlinear PDEs 
play  the  paramount role in establishing the
rate of convergence  of the approximations. The importance of
such estimates  is demonstrated recently  by \cite{Kr97},
\cite{Kr99} and \cite{Kr00}, presenting the first  rate of
convergence result in the sup norm of finite-difference 
approximations for fully nonlinear degenerate Bellman
equations.  Ideas from these publications are used and 
developed further in 
\cite{BJ02},  \cite{Ja03}, 
\cite{BJ05}, \cite{DK05a}, \cite{DK05b}, \cite{DK07}  and
\cite{GS08}. Recent results on estimating the Lipschitz 
constant and second order differences of finite-difference
approximations for a large class of fully nonlinear degenerate 
PDEs, including the normalized Bellman equations are presented
in 
\cite{Kr07}. In \cite{DK05b} first order derivatives of  
finite-difference 
approximations to 
{\em degenerate} linear parabolic and elliptic PDEs 
are estimated and are used 
 to establish sharp estimates on the rate of convergence 
of the approximations in the sup 
norm.

Finite-difference methods for solving PDEs have been extensively 
studied since the first half of the last century.  
Let us mention the pioneering papers by 
R. Courant, K.O. Friedrichs and 
H. Lewy \cite{CFL}, S. Gerschgorin \cite{Ge}, 
and publications  
by D.G. Aronson, J. Douglas,
F. John, H.O. Kreiss, O. Ladyzhenskaya, P.D. Lax,  
W. Littman,  
Lyusternik, J. von Neumann, 
I.G. Petrovskii,  A.A. Samarskii, G. Strang, 
A.N. Tikhonov, V. Thom\'ee,  O.B. Widlund 
and many others (see,
e.g.,  \cite{A}, \cite{Dou}, \cite{J}, \cite{Kre},  
\cite{Lax}, \cite{Li}, \cite{Liu}, 
\cite{MMMS}, \cite{NR}, \cite{Pet}, \cite{TS}, 
 \cite{W} 
 and the references there.) We refer also to the review paper 
\cite{Th}, handbook \cite{Th1990}, and well-known monographs 
and textbooks for more information on the subject 
(\cite{C}, \cite{FW}, \cite{GR}, \cite{LT}, \cite{MG}, 
\cite{RM}, \cite{S},  \cite{Str}). 

The paper is organized as follows. 
The main results, 
Theorems \ref{theorem 4.7.2} 
and \ref{theorem 07.9.25.2} are presented 
in Section \ref{section results}. 
Here we   
formulate also a result, 
Theorem \ref{theorem 1.08.02.08}, on  
accelerated finite difference schemes, 
which we will prove 
in the continuation of this paper by using 
Theorem \ref{theorem 4.7.2}. 
As we have pointed out above, the idea of the
proof of Theorem \ref{theorem 1.08.02.08} is based on a power
expansion of $u_{h}$ in $h$. This idea was already applied
by the authors to show how to accelerate other approximation schemes
(see, for instance, \cite{GK1})
and it seems to the authors that it was never used before
in the framework of finite difference schemes in the sup norm
for {\em degenerate\/} 
elliptic and parabolic equations  
although much effort was applied to developing 
this and other
methods of improved approximation for uniformly
nondegenerate equations in domains (see, for instance
\cite{Bo} and the references therein).
It is worth saying that, in contrast with \cite{Bo}
and many other papers dealing with the expansion,
we do not use any information from the theory of PDE
and, as a matter of fact, the existence of smooth
solutions for degenerate elliptic and parabolic equations
follows directly from our results. 
 We deduce   
Theorem \ref{theorem 07.9.25.2} from 
Theorem \ref{theorem 4.7.2}, and conclude 
Section \ref{section results} with verifying the rather 
delicate conditions of the 
main results, Assumptions 
\ref{assumption 11.22.11.06}  and \ref{assumption 07.10.16.1},
for a class of examples, given in 
  Remark \ref{remark
07.10.17.4} before the formulation of the theorems. 
  The  proof of Theorem \ref{theorem 4.7.2} is given  in
Section \ref{section proof}, and the final section, Section
\ref{section 4.9.1},   is devoted to 
further discussions of Assumptions 
\ref{assumption 11.22.11.06}  and \ref{assumption 07.10.16.1}.

\mysection{Formulation of the main results}               
                                     \label{section results}

We take  some numbers $h_{0},T\in(0,\infty)$ 
and for  each number $h\in(0,h_{0}]$ 
we consider the integral equation 
\begin{equation}                                          \label{equation}
u(t,x)=g_h(x)
+\int_0^t\big(L_hu(s,x)+f_h(s,x)\big)\,ds, 
\quad (t,x)\in H_T
\end{equation}
for $u$, where  $g=g_h=g_{h}(x)$ and $f=f_h=
f_{h}(s,x)$  are
given  real-valued Borel functions of 
$x\in\bR^d$ and  
$(s,x)\in H_T=[0,T]\times\bR^d$,  
respectively, and 
 $L=L_h$ is a linear operator 
defined by 
$$
L_h\varphi(t,x)=L_h^0\varphi(t,x)
-c(t,x)\varphi(x), 
$$
$$
L^0\varphi(t,x)=L_h^0\varphi(t,x)=\frac{1}{h}
\sum_{\lambda\in\Lambda_{1}}
q_{\lambda}(t,x )\delta_{h,\lambda}\varphi(x) 
+\sum_{\lambda\in\Lambda_{1}} 
p_{\lambda}(t,x )\delta_{h,\lambda}\varphi(x) , 
$$
for functions $\varphi$ on $\bR^d$. Here 
$\Lambda_{1}$ 
is a finite  subset of  $\bR^d$ 
such that
$0\not\in\Lambda_{1}$,   
and  $p_{\lambda}(t,x )$, $q_{\lambda}(t,x )$ are
real-valued functions  of 
$(t,x)\in H_{\infty} =[0,\infty)\times\bR^d $
given for each 
$\lambda\in\Lambda_{1}$, and
$$
\delta_{\lambda}\varphi(x )=\delta_{h,\lambda}\varphi(x )=
\frac{1}{h}(\varphi(x +h\lambda)-\varphi(x )),
\quad\lambda\in\Lambda_{1}.
$$  
As usual, for multi-indices $\alpha=(\alpha_1,\dots\alpha_d)$,
$\alpha_{i}=0,1,...$, we use the notation
$$
D^{\alpha}=D_{1}^{\alpha_{1}}...D_{d}^{\alpha_{d}},
\quad D_{i}=
\frac{\partial }{\partial x_i},\quad
|\alpha| =\sum_{i}\alpha_i, 
\quad
 D_{ij}=D_iD_j .
$$
For smooth $\varphi$ and integers $k\geq0$ we introduce
$D^{k}\varphi$ as the collection of partial derivatives
of $\varphi$ of order $k$, and define 
$$
|D^{k}\varphi|^{2}=\sum_{|\alpha|=k}|D^{\alpha}
\varphi|^{2},\quad
[\varphi]_{k}
=\sup_{x\in\bR^d}
|D^{k}\varphi(x)|,\quad |\varphi|_{k}=
\sum_{i\leq k}[\varphi]_{i}.
$$

Let $m\geq0$ be a fixed integer
and let $K_{1}\in[1,\infty)$ be
a constant. 
We make the following assumptions.

\begin{assumption}                  \label{assumption 16.12.07.06}
For any $\lambda\in\Lambda_{1}$ 
the derivatives in 
$x$ of  
  $p_{\lambda},q_{\lambda},c,f , g$ 
up to order $m$ 
are continuous functions in $(t,x)\in H_T$ 
and, for $k=0,...,m$ and some  
constants  $M_k$
we have
\begin{equation}
                                         \label{3.5.1}
\sup_{H_{T}} \big(
\sum_{\lambda\in\Lambda_{1}}
(|D^{k}q_{\lambda}|^{2}+|D^{k}p_{\lambda}|^{2}\big)
+|D^{k}c|^{2}\big)
\leq M^{2}_{k}. 
\end{equation} 
\end{assumption}

By Theorem 2.3 of  
\cite{GK} under Assumption \ref{assumption 16.12.07.06}
for each $h\in(0,h_{0}]$,
there exists a unique bounded solution $u_{h}$
of \eqref{equation} and this solution is
continuous in $H_{T}$
along with all its derivatives in $x$ up to order
$m$. 
However, the bounds, provided by this theorem 
for these  derivatives depend on the parameter $h$. 
Our aim is to show  the existence of bounds, independent 
of $h$, if in addition to  
Assumption \ref{assumption 16.12.07.06}, the assumptions below 
also hold.

\begin{assumption}            \label{assumption 18.12.07.06} 
For   all $t\in[0,T]$ 
\begin{equation}
                                         \label{3.9.1}
\sum_{\lambda\in\Lambda_{1}}\lambda q_{\lambda}(t,x )
\quad\hbox{is independent of}\quad  
 x .
\end{equation}
\end{assumption}
Introduce
$$
 \chi_{\lambda}= 
\chi_{h,\lambda}=q_{\lambda}+hp_{\lambda}.
$$
\begin{assumption}  
                            \label{assumption 1.26.11.06}   
For all $(t,x)\in H_T$, $h\in(0,h_{0}]$,
and $\lambda\in\Lambda_{1}$,
\begin{equation}                     \label{1.10.02.08}
\chi_{\lambda}(t,x )\geq 0.
\end{equation}
There exists a constant $c_{0}>0$ such that $c\geq c_{0}$.

\end{assumption}  

Obviously Assumption \ref{assumption 1.26.11.06}
implies that $q_{\lambda}\geq0$.

\begin{remark}
                                     \label{remark 07.9.18.8}
The above assumption: $c\geq c_{0}>0$, is almost irrelevant
if we only consider \eqref{equation} on a finite time interval.
Indeed,  if $c$ is just bounded, say $|c|\leq C=\text{const}$, by 
introducing a new function $v(t,x)=u(t,x)e^{-2Ct}$
we will have an equation for  $v$ similar to 
\eqref{equation} with $L^{0}v-(c+2C)v$ and $fe^{-2Ct}$
in place of $Lu$ and $f$, respectively. Now for the new $c$ we have
  $c+2C\geq C$.
\end{remark}

Take a  function $\tau_{\lambda}$ defined on 
$\Lambda_1$ taking values in $[0,\infty)$ 
 and 
for  $\lambda\in\Lambda_{1}$ 
introduce  
the operators
$$
T_{\lambda}\varphi=
T_{h,\lambda}\varphi(x)=\varphi(x+h\lambda),
\quad
\bar\delta_{\lambda}=\bar\delta_{h,\lambda}
=\tau_{\lambda}h^{-1}(T_{\lambda}-1).  
$$ 
For uniformity of notation we also introduce
$\Lambda_{2}$ as the set of fixed
 distinct vectors $\ell^{1},...,
\ell^{d}$ none of which is in $ 
\Lambda_{1}$ and   define
$$
\bar{\delta}_{\ell^{i}}=\bar{\delta} _{h,\ell^{i}}=\tau_{0}D_{i},
\quad 
 T_{\ell^i}= T_{h,\ell^{i}}=1,\quad
\Lambda=\Lambda_{1}\cup\Lambda_{2},
$$
where $\tau_{0}$ is a 
fixed parameter satisfying
$$ 
\tau_0>0.
$$
For integers $k=1,2,...$ and $\lambda^{i}\in\Lambda$,
$i=1,2,...,k$, introduce  the 
multi-vectors
$$
\lambda=(\lambda^{1},...,\lambda^{k})\in\Lambda^{k}
$$
and the operators
$$
 T_{\lambda}= T_{h,\lambda}=T_{h,\lambda^{1}}...
T_{h,\lambda^{k}},\quad
 \bar\delta_{\lambda}
=\bar\delta_{h,\lambda} 
= \bar\delta _{h,\lambda^{1} }...
 \bar\delta_{h,\lambda^{k}} .
$$
It is also
convenient to set
$\Lambda_{1}^{0}=\Lambda_{2}^{0}=\Lambda^{0}=\{0\}$ and
to introduce $\delta_{0}=\bar\delta_0$ and $T_{0}$ as unit
operators.
For $\mu\in\Lambda^{k}$ and $k\leq m$ we set
$$
Q\varphi =h^{-1}\sum_{\lambda\in\Lambda_{1}}
q_{\lambda} \delta_{\lambda}
\varphi,
\quad
Q_{\mu}\varphi =h^{-1}\sum_{\lambda\in\Lambda_{1}}
( \bar\delta _{\mu}q_{\lambda})\delta_{\lambda}
\varphi  ,
$$
$$
P\varphi = \sum_{\lambda\in\Lambda_{1}}
p_{\lambda} \delta_{\lambda}
\varphi,
\quad
P_{\mu}\varphi = \sum_{\lambda\in\Lambda_{1}}
 ( \bar\delta _{\mu}p_{\lambda} )\delta_{\lambda}
\varphi,
$$
$$
L^{0}_{ \mu} =Q_{ \mu} 
+P_{ \mu },
$$
$$
A_{k}(\varphi)=2\sum_{\lambda\in\Lambda^{k}}
( \bar\delta _{\lambda}
\varphi)L^{0}_{\lambda}T_{\lambda}\varphi,\quad
\cQ(\varphi)=\sum_{\mu\in\Lambda_{1} }
\chi_{\mu}
(\delta_{\mu}\varphi)^{2}.
$$

Below $B(\bR^{d})
$ is the set of  bounded Borel 
functions on
$\bR^{d}$ and $\frK$ is 
the set of bounded operators
$\cK=\cK_{h}=\cK_{h}(t)$ 
mapping $B(\bR^{d})
$ into itself preserving  the cone
of nonnegative functions 
and satisfying $\cK1\leq1$. 
 Set 
$$
|\Lambda_1|^2
=\sum_{\lambda\in\Lambda_1}|\lambda|^2,\quad
\|\Lambda_1\|^2
=\sum_{\lambda\in\Lambda_1}|\tau_{\lambda}\lambda|^2.
$$
 
Finally, fix a constant 
$\delta\in(0,1]$.

\begin{assumption}                \label{assumption 11.22.11.06} 
 
We have $m\geq1$ and 
for any $h\in(0,h_{0}]$, there
exists  an operator $\cK=\cK_{h,m}\in\frK$,
 such that  
\begin{equation}
                                             \label{3.24.1}
m A_{1}(\varphi)
\leq(1- \delta)\sum_{\lambda\in\Lambda}
\cQ( \bar\delta _{\lambda}\varphi)
+K_{1}\cQ(\varphi)
+2(1-\delta)c\cK\big(\sum_{\lambda\in\Lambda }
| \bar\delta _{\lambda}\varphi|^{2}\big) 
\end{equation}
on $H_{T}$
for all smooth functions $\varphi $.
\end{assumption}

\begin{assumption}         \label{assumption 07.10.16.1} 
We have $m\geq2$ and, for any $h\in(0,h_{0}]$
and $n=1,...,m$, there
exists  an operator $\cK=\cK_{h,n}\in\frK$,
 such that  
$$
n \sum_{\nu\in\Lambda}A_{1}( \bar\delta _{\nu}\varphi)
+ n(n-1)\sum_{\lambda\in\Lambda^{2}}
( \bar\delta _{\lambda}\varphi)
Q_{\lambda}T_{\lambda}\varphi
\leq (1- \delta)\sum_{\lambda\in\Lambda^{2}}
\cQ( \bar\delta _{\lambda}
\varphi)
$$
\begin{equation}
                                           \label{3.24.01}
+K_{1}
\sum_{\lambda\in\Lambda}
\cQ( \bar\delta _{\lambda}\varphi)
+2(1-\delta)c\cK\big(\sum_{\lambda\in\Lambda ^{2}}
| \bar\delta _{\lambda}\varphi|^{2}\big) 
+K_{1}\cK\big(\sum_{\lambda\in\Lambda}
| \bar\delta _{\lambda}\varphi|^{2}
\big)
\end{equation}
on $H_{T}$
for all smooth functions $\varphi $.
\end{assumption}
Obviously Assumptions 
\ref{assumption 11.22.11.06} and \ref{assumption 07.10.16.1} 
are satisfied if $q_{\lambda}$ and $p_{\lambda}$
are independent of $x$. In the general case,
as it is discussed in 
\cite{GK}, the above assumptions impose not only 
analytical conditions, but they are related also 
to some 
structural conditions, 
which can somewhat easier be analyzed under the following 
symmetry condition: 
\smallskip
 
 (S) $\Lambda_{1}=-\Lambda_{1}$ and
$q_{\lambda}=q_{-\lambda}$ 
for all $\lambda\in\Lambda_1$.
\smallskip
 
Notice that, if condition (S) holds then 
$$
h^{-1}\sum_{\lambda\in\Lambda_{1}}q_{\lambda}(t,x)
\delta_{\lambda}\varphi(x)
=(1/2)\sum_{\lambda\in\Lambda_{1}}q_{\lambda}(t,x)
\Delta_{\lambda}\varphi(x),
$$
where $\Delta_{\lambda}=\Delta_{h,\lambda}$ and
$$
\Delta_{h,\lambda}\varphi(x)
=\frac{\varphi(x+h\lambda)-2\varphi(x)+\varphi(x-h\lambda)}
{h^{2}}
=-\delta_{-\lambda}\delta_{\lambda}\varphi(x). 
$$
 
\begin{remark}                                     \label{remark 07.10.17.4} 
In this remark we
suppose that
Assumptions 
\ref{assumption 16.12.07.06}, 
\ref{assumption 18.12.07.06},
and \ref{assumption 1.26.11.06}  hold 
and $m\geq2$.
Assumption  \ref{assumption 11.22.11.06}
is discussed at length
and in many
details in \cite{GK}.
At the end of this section we show that 
if condition (S) holds  and ,
for all $\lambda\in\Lambda_{1}$, $\tau_{\lambda}>0$
and $q_{\lambda}\geq\kappa$,
where $\kappa>0$ is a constant, then
both Assumptions \ref{assumption 11.22.11.06}
and \ref{assumption 07.10.16.1} are satisfied 
for 
any $c_{0}>0$ 
 and $\delta\in(0,1)$,  
if
$h_{0}$ is sufficiently small and 
$ \tau_0  >0 ,K_{1}$,  and $\cK$
are chosen appropriately.
It follows from Remark  6.4 
of \cite{GK} and Remark 
\ref{remark 07.10.17.1} that the above  
 condition $\kappa>0$
can be dropped, provided,
additionally, that $c_{0}$
is large enough 
(this time we need not assume that $h$ is small).

 By the way, as we have seen in Remark \ref{remark 07.9.18.8},
the condition that $c_{0}$ be large is, actually,
harmless as long as we are concerned with equations
on a finite time interval.

Mixed situations, when $c$ is large
at those points where some of $q_{\lambda}$
can vanish are considered in Section
\ref{section 4.9.1}.
\end{remark} 

Now we are in the position to formulate our main results.

\begin{theorem}
                                          \label{theorem 4.7.2}
Let 
Assumptions \ref{assumption 16.12.07.06} 
through 
\ref{assumption 07.10.16.1} hold. 
Then for $h\in(0,h_{0}]$
we have 
\begin{equation}
                                             \label{4.8.03}
 \sup_{H_{T}}\sum_{k=0}^{m} 
|D^{k}u_h| \leq N 
 (F_{m}+G_{m}) ,
\end{equation}
where 
$$
 F_{n}
=\sum_{k\leq n} \sup_{H_{T}}
|D^{k}f_h| , 
\quad 
G_n=\sum_{k\leq n} \sup_{\bR^{d}}
|D^{k}g_h| ,
$$
  and $N$ 
depends only on   $\tau_{0}$,
$m$, $\delta$,  $c_0$, 
$K_{1}$, 
 $|\Lambda_1|$, $\|\Lambda_1\|$, $M_{0},...,M_{m}$. 
\end{theorem}

We prove this theorem in Section \ref{section proof}. Now we 
derive from it an estimate for the solution of the equation 

\begin{equation}                    \label{15.15.04.08}
L_hv+f_h=0\quad\text{in $\bR^d$},
\end{equation}
when $q_{\lambda}$, $p_{\lambda}$,
$c$, and $f$ are independent of $t$. 

\begin{theorem}
                                   \label{theorem 07.9.25.2}
Let  
Assumptions  \ref{assumption 16.12.07.06} 
through \ref{assumption 11.22.11.06}
 be satisfied. Suppose that $q_{\lambda}$, $p_{\lambda}$,
$c$, and $f$ are independent of $t$. Then the following statements 
hold. 

(i) There exists a unique bounded solution 
$v =v_h(x)$ of \eqref{15.15.04.08}. Moreover, all 
derivatives in $x$ of $v$ up to order $m$ are bounded continuous 
functions on $\bR^d$.

(ii) Let Assumption \ref{assumption 07.10.16.1} 
be also satisfied. Then  
\begin{equation}
                                             \label{4.8.03}
 \sup_{\bR^d}\sum_{k=0}^{m} 
|D^{k}v_h| \leq N F_{m},
\end{equation}
where 
$$
F_{m}=\sum_{k=0}^m\sup_{\bR^{d}}|D^kf_h|
$$
  and $N$
depends only on 
 $\tau_{0}$,
$m$, $\delta$,  $c_0$, 
$K_{1}$, 
 $|\Lambda_1|$, $\|\Lambda_1\|$, $M_{0},...,M_{m}$. 
\end{theorem}

\begin{proof} Statement (i) is proved in \cite{GK} (see 
Theorem 2.2 there). To prove (ii) take 
$\nu=c_{0}\gamma$,
where $\gamma>0$ is so small that  $c-\nu\geq c_{0}/2$ and conditions
\eqref{3.24.1} and \eqref{3.24.01}
 hold  with $c-\nu$ and $\delta/2$ in place
of $c$ and $\delta$, respectively.
Define 
$u(t,x):=v(x)e^{\nu t} $ and observe that
  $u$ satisfies 
$$
\frac{\partial}{\partial t}u
=L^{0}u-(c-\nu)u+e^{\nu t}f.
$$
By Theorem \ref{theorem 4.7.2} for $x\in \bR^d$ 
$$
e^{\nu T}\sum_{k=0}^m|D^kv(x)|
=\sum_{k=0}^m|D^ku(T,x)|
\leq
N e^{\nu T}F_{m}
+
N\sum_{k\leq m} \sup_{\bR^{d}}|D^kv(x)|.
$$
By multiplying the extreme terms by $e^{-\nu T}$
and letting $T\to\infty$, we get the result. 
\end{proof}
 
The above theorems have important applications in the numerical 
analysis of finite difference schemes for parabolic and elliptic 
PDEs. Using them in the continuation of the present 
paper we obtain accelerated finite difference schemes for 
second order (possibly) degenerate parabolic and also 
for second order (possibly) degenerate elliptic PDEs. In particular, 
we will consider the 
Cauchy problem  

\begin{equation}                          \label{1.07.02.08}
\frac{\partial}{\partial t}u(t,x)=
\cL u(t,x)+f_0(t,x), \quad t\in(0,T],\,x\in\bR^d
\end{equation}
\begin{equation}                          \label{2.07.02.08}
u(0,x)=g_0(x),\quad x\in\bR^d
\end{equation}
with the operator 
$$
\cL:=\tfrac{1}{2}
\sum_{\lambda\in\Lambda_1}
\sum_{i ,j=1}^d
q_{\lambda}\lambda_i\lambda_jD_iD_j
+\sum_{\lambda\in\Lambda_1}
\sum_{i=1}^dp_{\lambda}\lambda_iD_i-c.
$$
By a solution of \eqref{1.07.02.08}-\eqref{2.07.02.08}
we mean a continuous
 function $u(t,x)$ on $H_{T}$, such that for
each $t$ it is
  twice continuously differentiable in $x$, is
bounded in $H_{T}$ along with 
 its  derivatives 
 in $x$ up to second order 
and satisfies
$$
u(t,x)=g_{0}(x)+\int_{0}^{t}[\cL u(s,x)+
f_{0}(s,x)]\,ds 
$$
in $H_{T}$. 

To formulate one of the main results of the continuation 
of the paper  we fix an integer $k\geq0$
and  set $f_{h}=f_{0}$,  $g_{h}=g_{0}$, 
$$
\bar u_h=\sum_{j=0}^{k}b_ju_{2^{-j}h}, 
\quad 
$$
where $u_{2^{-j}h}$ is 
the solution to \eqref{equation} 
with $2^{-j}h$ in place of $h$, 
$$
(b_0,b_1,...,b_k)
:=(1,0,0,...,0)V^{-1}, 
$$
and $V^{-1}$ is the inverse of the 
Vandermonde matrix with entries
$$
V^{ij}:=2^{-(i-1)(j-1)}, \quad i,j=1,...,k+1.
$$

\begin{theorem}                  \label{theorem 1.08.02.08}
Let Assumptions 
\ref{assumption 1.26.11.06} and 
\ref{assumption 16.12.07.06}  
with $m\geq 3(k+1)$
hold. Also let  
  condition (S) be satisfied.  
Then \eqref{1.07.02.08}--\eqref{2.07.02.08} has a 
unique solution $u_0$, and 
\begin{equation}                    \label{2.08.02.08}                
|\bar u_h(t,x)-u_0(t,x)|\leq N h^{k+1}
\end{equation}                  
holds for all $(t,x)\in H_T$, 
$h\in(0,h_0]$, where 
$N$ is a constant depending only on $T$,
$k$, $d$, $|\Lambda_1|$, $h_0$, 
the number of elements 
in $\Lambda_1$, on $M_0,\dots. M_m$, on 
$\sup_{t\in[0,T]}|f_0(t)|_m$ and on $|g_0|_m$. 
\end{theorem}

Now we prove Remark \ref{remark 07.10.17.4}.  
Instead of condition (S) we assume the following
weaker condition 
\smallskip

(S$'$): 
$\Lambda_1=-\Lambda_1$ and $Dq_{\lambda}=Dq_{-\lambda}$ 
for $\lambda\in\Lambda_1$, 
\smallskip

\noindent
and proceed with the proof of Remark \ref{remark 07.10.17.4} 
as follows.  
Clearly, 
$$
\sum_{\lambda\in\Lambda}(
\bar{\delta}_{\lambda}\varphi)
L^0_{\lambda}T_{\lambda}\varphi=
I_1+I_2, 
$$
with 
$$
I_1:=\sum_{\lambda\in\Lambda_1}(\bar{\delta}_{\lambda}\varphi)
L^0_{\lambda}T_{\lambda}\varphi, 
\quad
I_2:=\sum_{\lambda\in\Lambda_2}(\bar{\delta}_{\lambda}\varphi)
L^0_{\lambda}\varphi.
$$
Due to condition (S$'$)
$$ 
I_1=
\sum_{\lambda\in\Lambda_1}(\bar\delta_{\lambda}\varphi)
L^0_{\lambda}\varphi
+
h\sum_{\lambda\in\Lambda_1}
(\bar\delta_{\lambda}\varphi)
L^0_{\lambda}\delta_{\lambda}\varphi
$$
$$
=\tfrac{1}{2}\sum_{\lambda ,\mu\in\Lambda_1}
(\bar\delta_{\lambda}\varphi)
(\bar\delta_{\lambda}q_{\mu})\Delta_{\mu}\varphi
+\sum_{\lambda ,\mu\in\Lambda_1}
(\bar\delta_{\lambda}\varphi)
(\bar\delta_{\lambda}p_{\mu})\delta_{\mu}\varphi
$$
$$
+\sum_{\lambda ,\mu\in\Lambda_1}
(\bar\delta_{\lambda}\varphi)
(\bar\delta_{\lambda}\chi_{\mu})\delta_{\mu}
\delta_{\lambda}\varphi
=:I_1^{(1)}+I_1^{(2)}+I_1^{(3)}, 
$$
$$
I_2=I_2^{(1)}+I_2^{(2)},
$$
where in the notation $\xi=D\varphi/|D\varphi|$ 
and 
$\psi_{(\xi)}=\xi_{i}D_{i}\psi$,
$$
I_2^{(1)}=
\tfrac{1}{2}\tau_{0}^{2}\sum_{j=1}^{d}\sum_{ \mu\in\Lambda_1}
(D_{j}\varphi)
(D_{j}q_{\mu})\Delta_{\mu}\varphi
=\tfrac{1}{2}\tau_{0}^{2}|D\varphi|
\sum_{\mu\in\Lambda_{1}}q_{\mu(\xi)}\Delta_{\mu}\varphi,
$$
$$
I_2^{(2)}=
 \tau_{0}^{2}\sum_{j=1}^{d}
\sum_{ \mu\in\Lambda_1}
(D_{j}\varphi)
(D_{j}p_{\mu})\delta_{\mu}\varphi=
\tau_{0}^{2}|D\varphi|
\sum_{ \mu\in\Lambda_1}
 p_{\mu(\xi)} \delta_{\mu}\varphi.
$$
Set
$$
\bar\tau:=\max_{\lambda\in\Lambda_1}\tau_{\lambda},
\quad
\underline\tau:=\min_{\lambda\in\Lambda_1}\tau_{\lambda},
$$ 
and observe that $\chi_{\lambda}\geq\kappa/2>0$ for sufficiently small 
$h$, and that $c\geq c_0>0$. Then
by Young's inequality,  
we obtain 
$$
2mI^{(j)}_1
\leq\tfrac{1-\delta}{3} \sum_{\lambda\in\Lambda_1}
\cQ(\bar\delta_{\lambda} \varphi)
+N\cQ(\varphi)\quad \text{for $j=1,3$}, 
\quad 2mI^{(2)}_{1}\leq N\cQ(\varphi),
$$
$$
2mI^{(1)}_2\leq \tfrac{1-\delta}{3}\sum_{
\mu\in\Lambda_1}
\cQ(\bar\delta_{\mu}\varphi)
+ \tau_{0}^{2}Nc_0^{-1}c
\sum_{\lambda\in\Lambda_2}|\bar\delta_{\lambda}\varphi|^2, 
$$ 
$$
2mI^{(2)}_2\leq \tau_{0}Nc_0^{-1}c
\sum_{\lambda\in\Lambda}|\bar\delta_{\lambda}\varphi|^2,
$$
where 
$N$ is a constant depending only on $m$, 
$\kappa$, $\delta$, 
$\bar\tau$,
$\underline\tau$, 
the number of elements in $\Lambda_1$,
and on the supremum norm of 
the gradients of 
$p_{\lambda}$ and $q_{\lambda}$ in $x$. 
Summing up these inequalities 
and taking $\tau_{0}>0$ sufficiently 
small  we get 
\eqref{3.24.1} with $K_1=3N$,
  unit operator $\cK$, and with 
$\delta$ as close to $1$ as we wish. 

This result is obviously applicable
to $\bar{\delta}_{\nu}\varphi$ is place of $\varphi$
for any $\nu\in\Lambda$. It follows that
for any $\delta\in(0,1)$  and appropriate 
constant $ K_1 $
we have
$$
m\sum_{\nu\in\Lambda}A_{1}( \bar\delta _{\nu}\varphi)
\leq(1- \delta)\sum_{\lambda\in\Lambda^{2}}
\cQ( \bar\delta _{\lambda}\varphi)
$$
\begin{equation}
                                        \label{08.4.27.1}
+K_{1}\sum_{\lambda\in\Lambda}\cQ(\bar\delta _{\lambda}\varphi)
+2(1-\delta)c \sum_{\lambda\in\Lambda^{2} }
| \bar\delta _{\lambda}\varphi|^{2} .
\end{equation}

Now we show that 
Assumption 
\ref{assumption 07.10.16.1} holds. 
Clearly, 
$$
\sum_{\lambda\in\Lambda^{2}}
(\bar\delta_{\lambda}\varphi)
Q_{\lambda}T_{\lambda}\varphi=
\sum_{\lambda\in\Lambda_1^{2}}
(\bar\delta_{\lambda}\varphi)
Q_{\lambda}T_{\lambda}\varphi+
2\sum_{\lambda\in\Lambda_1\times\Lambda_2}
(\bar\delta_{\lambda}\varphi)
Q_{\lambda}T_{\lambda}\varphi
$$
$$
+\sum_{\lambda\in\Lambda_2^{2}}
(\bar\delta_{\lambda}\varphi)
Q_{\lambda}\varphi=I_1+I_2+I_3. 
$$
Using  
$$
\Delta_{\mu}T_{\lambda}=\Delta_{\mu}
+h\delta_{\lambda^{1}}\delta_{\lambda^{2}}(\delta_{\mu}
+\delta_{-\mu})+(
\delta_{\lambda^{1}}
+\delta_{\lambda^{2}})(\delta_{\mu}
+\delta_{-\mu})
$$
for $\lambda=(\lambda^{1},\lambda^{2})\in\Lambda_1^2$ and
$\mu\in\Lambda_1$,  we have  
$$
I_1=\sum_{\lambda\in\Lambda_{1}^{2}}
(\bar\delta_{\lambda}\varphi)
Q_{\lambda}T_{\lambda}\varphi
= I^{(1)}_{1}+I_1^{(2)},
$$
with 
$$
I_{1}^{(1)}=
\tfrac{1}{2}
\sum_{\lambda\in\Lambda_{1}^{2}\mu\in\Lambda_1}
(\bar\delta_{\lambda}\varphi)
(\bar\delta_{\lambda}q_{\mu}) (4
\delta_{\lambda^{1}}-\delta_{-\mu})\delta_{\mu}\varphi,
$$
$$
I_1^{(2)}=
h\sum_{\lambda\in\Lambda_{1}^{2},\mu\in\Lambda_1}
(\bar\delta_{\lambda}\varphi)
(\bar\delta_{\lambda}q_{\mu})
\delta_{\lambda} \delta_{\mu}
\varphi.
$$ 
As above,
we have 
$$
\sum_{\lambda\in\Lambda_1^2}
|\bar\delta_{\lambda}\varphi|^2
\leq N
\sum_{\lambda,\mu\in\Lambda_1}
|\delta_{\lambda}\delta_{\mu}\varphi|^2
\leq N
\sum_{\lambda\in\Lambda}
\cQ(\bar\delta_{\lambda}\varphi), 
$$
$$
\sum_{\lambda\in\Lambda_1^2,\mu\in\Lambda_1}
|\delta_{\lambda}\delta_{\mu}\varphi|^2
\leq N
\sum_{\lambda\in\Lambda^2}
\cQ(\bar\delta_{ \lambda }\varphi), 
$$
and hence, by Young's inequality
\begin{equation}                             \label{1.05.02.08}
n(n-1)I_{1}^{(1)}
\leq N \sum_{\lambda\in\Lambda}
\cQ(\bar\delta_{\tau}\varphi), 
\end{equation}
\begin{equation}                             \label{2.05.02.08}
n(n-1)I_{1}^{(2)}
\leq h^2N \sum_{\lambda\in\Lambda^2}
\cQ(\bar\delta_{\lambda}\varphi)
 +N \sum_{\lambda\in\Lambda}
\cQ(\bar\delta_{\lambda}\varphi), 
\end{equation}
where $N$ is a constant depending only on 
$m$, 
$\kappa$, 
$\bar\tau$, 
the number of elements in $\Lambda_1$ 
and on the supremum norm of 
$|D^2q_{\lambda}|$. 
Similarly,
$$             
n(n-1)I_2=n(n-1)\tau_0\sum_{i=1}^d
\sum_{\nu,\mu\in\Lambda_1}
(\bar\delta_{\nu}\tau_0D_i\varphi)
(\bar\delta_{\nu}D_i q_{\mu})
T_{\nu}\Delta_{\mu}\varphi
$$
\begin{equation}                            \label{3.05.02.08}  
\leq 
\tau_0N\sum_{\lambda\in\Lambda}
\cQ(\bar\delta\varphi)
+N\tau_0c
\cK(\sum_{\mu\in\Lambda^2}
\bar\delta_{\mu}\varphi), 
\end{equation}
$$  
n(n-1)I_3=\tfrac{1}{2}n(n-1)
\tau_0^2\sum_{i,j=1}^d\sum_{\mu\in\Lambda_1}
(\tau^2_0D_{ij}\varphi)
(D_{ij}q_{\mu})\Delta_{\mu}\varphi
$$
\begin{equation}                            \label{4.05.02.08}
\leq N\tau_0^2c
\sum_{\lambda\in\Lambda^2}
|\bar\delta_{\lambda}\varphi|^2,
\end{equation}
where $N$  denote  some  constants depending on 
$m$, $d$, $c_0$, 
$\kappa$, 
$\bar\tau$, $\underline\tau$,
the number of elements in $\Lambda_1$ 
and on the supremum norm of 
$|D^2q_{\lambda}|$. 
Summing up the inequalities 
\eqref{08.4.27.1} through  
\eqref{4.05.02.08} 
and choosing $\tau_0$ and $h_0$ sufficiently 
small we obtain \eqref{3.24.01}.

\mysection{Proof of Theorem \protect\ref{theorem 4.7.2}}    
                                         \label{section proof}
 
For $m=1$ estimate \eqref{4.8.03} holds by 
virtue of Theorem 2.1 from 
\cite{GK}, proved 
by the aid of the following version 
of the maximum principle   
(Corollary 3.2 in \cite{GK}).

\begin{lemma}
                                        \label{lemma 3.6.1}
Let Assumption
\ref{assumption 16.12.07.06} with $m=0$ be satisfied and let
$\chi_{\lambda}\geq0$ for all $\lambda\in\Lambda_{1}$.
Let $v$ be a bounded
 function on $H_{T}$, such that 
the partial derivative $D_{t}v:=\partial v(t,x)/\partial t$
exists in $H_{T}$. 
Let   
$F$ be a nonnegative integrable function on $[0,T]$,
and let $C$ be a nonnegative bounded function on $H_T$ 
such that 
$$
\nu:=\sup_{H_{T}}(C-c)<0. 
$$ 
Assume that
  for all $(t,x)\in H_{T}$ we have
\begin{equation}
                                               \label{07.9.23.1}
D_{t}v \leq  L v + C\bar{v}_{+}+F,
\end{equation}
where $\bar{v} (t)=\sup\{v (t,x):x\in \bR^d\}$.
Then in $[0,T]$ we have
\begin{equation}
                                        \label{3.11.1}
\bar{v}(t)\leq  
\bar v_{+}(0)
+|\nu|^{-1}\sup_{[0,t]}F,
\end{equation}
where $a_+:=(|a|+a)/2$ for real numbers $a$. 
\end{lemma}

For the proof of this lemma we refer to \cite{GK}. 
In order to obtain Theorem \ref{theorem 4.7.2} for $m\geq2$ 
we need some 
more lemmas. 
First we prove 
 a lemma which will be used a few times 
in the future.
By $\cK$ in the lemma and later in the article
we mean a generic operator of class $\frK$. This operator
may change each time it is mentioned even in one line
(cf. the use of $o(n)$). 
  Thus, for example, for 
nonnegative 
functions $\alpha$, $\beta$ on $\bR^d$ 
the 
formula $\alpha\cK+\beta\cK=(\alpha+\beta)\cK$ 
means the 
simple fact that for any 
$\cK_1,\cK_2\in\frK$   
$$
\alpha\cK_1+\beta\cK_2=(\alpha+\beta)\cK_3
$$ 
with 
$$
\cK_3:=\tfrac{\alpha}{\alpha+\beta}\cK_1
+\tfrac{\beta}{\alpha+\beta}\cK_2
\in\frK\quad \big(\tfrac{0}{0}:=0\big).
$$

\begin{lemma}
                                   \label{lemma 3.21.3}
Let 
Assumption \ref{assumption 16.12.07.06} 
be satisfied.
Let $n\geq1$  be an  integer. 
Set 
$$
A^{2}=\|\Lambda_{1}\|^{2}+\tau_{0}^{2}.
$$ 
Then   for any
 $\varphi\in C^{n}$ 
we have 
\begin{equation}
                                            \label{3.21.5}
\sum_{\lambda\in\Lambda^{n}_{1}}
| \bar\delta_{\lambda}\varphi|^{2}
\leq
\|\Lambda_1\|^{2n} 
\cK\big( |D^{n}\varphi|^{2}\big),
\quad
\sum_{\lambda\in\Lambda^{n}}
| \bar\delta _{\lambda}\varphi|^{2}
\leq A^{2n} 
\cK\big( 
|D^{n}\varphi|^{2}\big).
\end{equation}
Furthermore, if $1\leq n\leq m$, then  
for any $\varphi\in C^{n}$
$$
\sum_{ 
\lambda\in\Lambda^{n}} |
P _{\lambda}T_{\lambda} \varphi|^{2}
\leq  |\Lambda_1|^2A^{2n}
\big(\sup_{H_T}\sum_{\mu\in\Lambda_1}|D^np_{\mu}|^2\big) 
\cK(|D\varphi|^{2})
$$
\begin{equation}
                                             \label{3.22.2}
\leq \tau_{0}^{-2}  |\Lambda_1 |^2A^{2n}
\big(\sup_{H_T}\sum_{\mu\in\Lambda_1}|D^np_{\mu}|^2\big) 
\cK\big(\sum_{\lambda\in\Lambda}|\bar{\delta}_{\lambda}
\varphi|^{2}\big),
\end{equation}
 and if assumption \eqref{3.9.1} holds,  
$$
\sum_{ 
\lambda\in\Lambda^{n}} |
Q_{\lambda}T_{\lambda}  \varphi|^{2} 
  \leq 
 |\Lambda_1|^4A^{2n}
\big(\sup_{H_T}\sum_{\mu\in\Lambda_1}|D^nq_{\mu}|^2 \big)\cK(  
|D^{2}\varphi|^{2})
$$
\begin{equation}
                                             \label{3.22.3}  
 \leq \tau_{0}^{-4}
 |\Lambda_1|^4 A^{2n}\big
(\sup_{H_T}\sum_{\mu\in\Lambda_1}|D^nq_{\mu}|^2\big )\cK\big(  
\sum_{\lambda\in\Lambda^{2}}|\bar{\delta}_{\lambda}\varphi|^{2}).
\end{equation} 
 Finally,
$$
\big(\sum_{\lambda\in\Lambda_{1}}
q_{\lambda}(\delta_{\lambda}\varphi)^{2}\big)^{2}  
\leq M_{0}^{2}
\sum_{\lambda\in\Lambda_{1}}
 (\delta_{\lambda}\varphi)^{4}
\leq M_{0}^{2}\big(
\sum_{\lambda\in\Lambda_{1}}
 (\delta_{\lambda}\varphi)^{2}\big)^{2} .
$$
\end{lemma} 
\begin{proof}
It is easy to see that
for $\lambda\in\Lambda^{n}_{1}$ we have
$$
 \bar\delta _{\lambda}\varphi(x)
=h^{-n}\int_{[0,h]^{n}}
\varphi_{\lambda }(\theta,x)\,d\theta,
$$
where, for $y(\lambda,\theta)= \lambda^{1}
\theta^{1}+...+\lambda^{n}\theta^{n}$ 
 and $\tau_{\lambda}
=\tau_{\lambda^1}\tau_{\lambda^2}
\cdot...\cdot \tau_{\lambda^n}$, 
$$
\varphi_{\lambda }(\theta,x)=\tau_{\lambda} 
\sum_{i_{1},...,i_{n}=1}^{d} 
\lambda^{1}_{i_{1}}
\cdot
...
\cdot
\lambda^{n}_{i_{n}}D_{i_{1}}\cdot...
\cdot D_{i_{n}}\varphi(x+y(\lambda,\theta)).
$$
By Cauchy's inequality 
$$
| \bar\delta _{\lambda}\varphi(x)|^{2}\leq
h^{-n}\int_{[0,h]^{n}}
|\varphi_{\lambda }(\theta,x)|^{2}\,d\theta,
$$
$$
|\varphi_{\lambda }(\theta,x)|\leq
 \tau_{\lambda} |\lambda | \,
|D^{n}\varphi(x+y(\lambda,\theta))|,
$$
where $|\lambda| :=|\lambda^{1}| \cdot...\cdot
|\lambda^{n}| $. It follows that the first inequality
in \eqref{3.21.5} holds with with $\cD_{h}$
in place of $\cK$, where 
$$
\cD_{h}\psi(x)=
\|\Lambda_1\|^{-2n} 
\sum_{\lambda\in\Lambda_{1}^{n}}
 \tau_{\lambda}^2 |\lambda|^{2}
h^{-n}\int_{[0,h]^{n}}
\psi(x+y(\lambda,\theta))\,d\theta
\quad  \big(\tfrac{1}{0}:=0\big).
$$
Since 
$$
\sum_{\lambda\in\Lambda_{1}^{n}}
\tau_{\lambda}^2|\lambda|^{2}=\|\Lambda_1\|^{2n}, 
$$ 
$
\cD_{h}1\leq1,
$  
that is $\cD_{h}\in\frK$ and  
the first inequality
in \eqref{3.21.5} is proved. 

To prove the second one 
introduce $\cD_{h,k}$ as the operators for which
first inequality
in \eqref{3.21.5} holds with $k$ in place of $n$,
recall that
$\Lambda_1^0=\Lambda_2^0=\{0\}$ and
$\bar\delta_0$ is the identity operator, 
and
observe that the left-hand side of the second
inequality equals
$$
\sum_{k= 0 }^{n}C^{k}_{n}
\sum_{\lambda\in\Lambda_{1}^{k},\mu\in
\Lambda_{2}^{n-k} }
| \bar\delta _{\lambda}
 \bar\delta _{\mu}
\varphi|^{2}\leq
\sum_{k= 0 }^{n}C^{k}_{n}
\|\Lambda_1\|^{2k} 
\cD_{h,k}
\big(
\sum_{ \mu\in
\Lambda_{2}^{n-k}}
 | \bar\delta _{\mu}D^{k}\varphi|^{2}
\big)
$$
$$
=\sum_{k= 0 }^{n}C^{k}_{n}
\|\Lambda_1\|^{2k}\tau_0^{2(n-k)} \cD_{h,k}
\big(
 | D^{n}\varphi|^{2}\big)=:
 (\|\Lambda_1\|^{2}+\tau_0^2)^n 
\cE_{h}\big(
 | D^{n}\varphi|^{2}\big),
$$
 with  $\cE_{h}\in\frK$.
This proves the first 
assertion of the lemma.

To prove \eqref{3.22.2} notice that
by  Cauchy's inequality and \eqref{3.21.5}
for $ \lambda\in\Lambda^{n}$
$$
|P _{\lambda}T_{\lambda}  \varphi|^{2}
=
\big|\sum_{\mu\in\Lambda_{1}}
( \bar\delta _{\lambda}
p_{\mu})T_{ \lambda} 
\delta_{\mu} \varphi\big|^{2}
$$
 $$
\leq \sum_{\mu\in\Lambda_{1}}
( \bar\delta _{\lambda}p_{\mu})^{2}
T_{\lambda} \sum_{\mu\in\Lambda_{1}}
(
\delta_{\mu} \varphi)^{2}
\leq  |\Lambda_1|^2 
\sum_{\mu\in\Lambda_{1}}
( \bar\delta _{\lambda}p_{\mu})^{2}
T_{\lambda} \cD_{h,1}(|D\varphi|^{2}).
$$
Hence the left-hand side 
of \eqref{3.22.2} is less
than
$$
 |\Lambda_1|^2 
\sum_{ \mu\in\Lambda_{1},
\lambda\in\Lambda^{n}}
( \bar\delta _{\lambda}p_{\mu})^{2}
T_{\lambda} \cD_{h, 1 }(|D\varphi|^{2})
$$
$$
=:
|\Lambda_1|^2A^{2n}
\big(\sup_{H_T}\sum_{\mu\in\Lambda_1}|D^np_{\mu}|^2\big)
\cH_{h}(|D\varphi|^{2}).
$$
Here $\cH_{h}\in\frK$,  since by \eqref{3.21.5}
$$
\sum_{ \mu\in\Lambda_{1},
\lambda\in\Lambda^{n}}
( \bar\delta _{\lambda}p_{\mu})^{2}\leq
A^{2n} 
\cD_{h,n}\big(\sum_{ \mu\in\Lambda_{1}}|D^{n}
p_{\mu}|^{2}\big)\leq 
A^{2n}
\sup_{H_T}\sum_{\mu\in\Lambda_1}|D^np_{\mu}|^2. 
$$
This proves   \eqref{3.22.2}.
To prove \eqref{3.22.3},
notice that 
  $\sum_{\mu}\mu \delta_{\lambda}
q_{\mu}=0$, which implies that
$$
Q_{\lambda}T_{\lambda}\varphi
=h^{-1}\sum_{\mu\in\Lambda_{1}}
( \bar\delta _{\lambda}q_{\mu}
)T_{\lambda}
(\delta_{\mu}\varphi-\mu_{i}D_{i}\varphi)
=\sum_{\mu\in\Lambda_{1}}
( \bar\delta _{\lambda}q_{\mu}
)T_{\lambda}\psi_{ \mu},
$$
where
$$
\psi_{ \mu}=h^{-1}
(\delta_{\mu}\varphi-\mu_{i}D_{i}\varphi).
$$
Hence as above
the left-hand side of \eqref{3.22.3} is less than
$$
\sum_{\lambda\in\Lambda^{n}}
\sum_{\mu\in\Lambda_{1}}
( \bar\delta _{\lambda}q_{\mu}
)^{2}T_{\lambda}\sum_{\mu\in\Lambda_{1}}|\psi_{\mu}|^{2}
=:
A^{2n}
\big(\sup_{H_T}\sum_{\mu\in\Lambda_1}|D^nq_{\mu}|^2
\big) \cF_{h}\big(
\sum_{\mu\in\Lambda_{1}}|\psi_{\mu}|^{2}\big),
$$
where $\cF_{h}\in\frK$. Furthermore,
$$
\psi_{\mu}
(x)=h^{-2}\int_{0}^{h}(h-\theta)\mu_{i}\mu_{j}
D_{ij}\varphi(x+\mu\theta)\,d\theta
$$
$$
\leq|\mu|^{2}h^{-1}
\int_{0}^{h}|
D^{2}\varphi(x+\mu\theta)|\,d\theta
 =|\mu|^{2}
\int_{0}^{1}|
D^{2}\varphi(x+h\mu\theta)|\,d\theta ,
$$
$$
|\psi_{\mu}
(x)|^{2}\leq|\mu|^{4}
\int_{0}^{1}|
D^{2}\varphi(x+h\mu\theta)|^{2}\,d\theta,
$$
and we obtain \eqref{3.22.3} with the operator
$$
\cK\psi:= |\Lambda_1|^{-4} 
\cF_{h}\big(
\sum_{\mu\in\Lambda_{1}}
|\mu|^{4}
\int_{0}^{1}\psi(\cdot+h\mu\theta) \,d\theta\big),
$$
which is in $\frK$ because
$$
\cK 1= |\Lambda_1|^{-4} 
\sum_{\mu\in\Lambda_{1}}
|\mu|^{4}\leq1.
$$
Since the last assertion of the lemma
is obvious, the lemma is proved. 
\end{proof}

The following lemma can be proved easily by induction
on $n$. 
 (Sums over empty sets of indices are defined 
to be 0 in the lemma, and everywhere in the article.)
 
\begin{lemma}
                                   \label{lemma 4.7.2}
Let $n\geq1$ be an integer, $\psi $ 
  and 
$\varphi$
be $n$ times continuously
differentiable functions
on $\bR^{d}$, and $\lambda\in\Lambda^{n}$. 
Then
$$
 \bar\delta _{\lambda}(\psi\varphi)
=\psi \bar\delta _{\lambda}\varphi
+\sum_{i=1}^{n} (\bar\delta _{\lambda^{i}}\psi)
 \bar\delta _{\bar\lambda (i)}T_{\lambda^{i}}\varphi
$$
$$
+\sum_{1\leq i<j\leq n} ( 
 \bar\delta _{\lambda(i,j)}\psi)
 \bar\delta _{\bar{\lambda}
(i,j)} T_{\lambda(i,j)}
\varphi+.... 
$$
$$
+\sum_{1\leq i_{1}<...<i_{k}\leq n} 
( \bar\delta _{\lambda(i_{1},...,i_{k})}\psi)
 \bar\delta _{\bar{\lambda}(i_{1},...,i_{k})}
T_{\lambda(i_{1},...,i_{k})}
\varphi
$$
\begin{equation}
                                            \label{4.7.3}
+...+( \bar\delta _{\lambda}\psi)
T_{\lambda}\varphi,
\end{equation}
where $\lambda(i_{1},...,i_{k})=(\lambda^{i_{1}},
....,\lambda^{i_{k}})$, $\bar{
\lambda}(i_{1},...,i_{k})$ is the sequence
of vectors $\lambda^{1},...,\lambda^{n}$ from which
the vectors standing on the places with
numbers $ i_{1} ,..., i_{k} $ 
are removed,  
 and $\bar\delta_{\bar\lambda(1) }:=1$ 
for $n=1$. 

\end{lemma}

{\bf Proof of Theorem \ref{theorem 4.7.2}}. 
Recall that 
$\Lambda^{0}=\{0\}$ and  
$ \bar\delta _{0}=T_{0}$
is the {\em unit\/} operator. 
Fix $h\in(0,h_{0}]$, for $0\leq k\leq m$  set
$$
u=u_{h},\quad f=f_{h},\quad
V_0=u^2, 
\quad
V_{k }= \sum_{\lambda
\in\Lambda^{k}}| \bar\delta _{\lambda}u|^{2},
\quad 
\bar{V}_{k }(t)=\sup_{\bR^{d}}V_{k }(t,x),
$$
and recall that $F_{n}$ is introduced
in Theorem \ref{theorem 4.7.2}. 
Take an integer $n\in[1,m]$.
Then we have
\begin{equation}                       \label{20.23.01.07}
L^{0}_{h}V_{n}=2 \sum_{\lambda\in\Lambda^{n}}
( \bar\delta _{\lambda}u)L^{0}_{h}
 \bar\delta_{\lambda} u
+\sum_{\lambda\in\Lambda^{n}}
\cQ( \bar\delta _{\lambda}u).
\end{equation}
By Lemma \ref{lemma 4.7.2}
\begin{equation}
                                         \label{4.8.3}
2 \sum_{\lambda\in\Lambda^{n}}
( \bar\delta _{\lambda}u)
L^{0}_{h} \bar\delta _{\lambda}u  
=2 \sum_{\lambda\in\Lambda^{n}}
( \bar\delta _{\lambda}u)
 \bar\delta _{\lambda}L^{0}_{h}u
-\sum_{n\geq k\geq1}I_{ n,k},
\end{equation}
where 
$$
I_{n,k}:= C^{k}_{n}
\sum_{\mu\in\Lambda^{n-k}}A_{k}( \bar\delta _{\mu}u).
$$
By Assumption \ref{assumption 11.22.11.06},
$$
nA_{1}(\bar\delta _{\mu}u)
\leq(1- \delta)
\sum_{\lambda\in\Lambda}
\cQ( \bar\delta _{\lambda}
 \bar\delta _{\mu} u)
+K_{1}\cQ( \bar\delta _{\mu}u)
$$
\begin{equation}                           \label{21.23.01.07}
+2(1-\delta)c\cK
\big(\sum_{\lambda\in\Lambda}
| \bar\delta _{\lambda} \bar\delta _{\mu}
u|^{2}\big).
\end{equation} 
Hence,
$$
I_{n,1}\leq (1- \delta)\sum_{\lambda\in\Lambda^{n}}
\cQ( \bar\delta _{\lambda}  u)
+K_{1}\sum_{\lambda\in\Lambda^{n-1}}
\cQ( \bar\delta _{\lambda}  u)+ 2(1-\delta)c\bar{V}_{n}.
$$
Next, if $n\geq2$, then
$$
I_{n,1}+I_{n,2}=
n\sum_{\mu\in\Lambda^{n-1}}A_{1}
( \bar\delta _{\mu}u)
+ \tfrac{1}{2} n(n-1)\sum_{\mu\in\Lambda^{n-2}}A_{2}
( \bar\delta _{\mu}u)
$$
$$
=\sum_{\mu\in\Lambda^{n-2}}
\big(n\sum_{\nu\in\Lambda}A_{1}
( \bar\delta _{\nu}( \bar\delta _{\mu}u))
+ \tfrac{1}{2} 
n(n-1)A_{2}( \bar\delta _{\mu}u)\big),
$$
so that
by Assumption \ref{assumption 07.10.16.1}
$$
I_{n,1}+I_{n,2}\leq  n(n-1)\sum_{\mu\in\Lambda^{n-2}}
\sum_{\lambda\in\Lambda^{2}}
( \bar\delta _{\lambda}
 \bar\delta _{\mu}u)
P_{\lambda}T_{\lambda}
 \bar\delta _{\mu}u
$$
$$
+(1- \delta)\sum_{\lambda\in\Lambda^{n}}
\cQ( \bar\delta _{\lambda} u)
+K_{1}
\sum_{\lambda\in\Lambda^{n-1}}
\cQ( \bar\delta _{\lambda} u)
+2(1-\delta)c\bar{V}_{n}+K_{1}\bar{V}_{n-1}.
$$
By Lemma \ref{lemma 3.21.3} 
$$
 n(n-1)\sum_{\mu\in\Lambda^{n-2}}
\sum_{\lambda\in\Lambda^{2}}
( \bar\delta _{\lambda}
 \bar\delta _{\mu}u)P_{\lambda}T_{\lambda}
 \bar\delta _{\mu}u+
\sum_{n\geq k\geq3}I_{ n, k}
\leq  \delta^{2} 
V_{n}+N \sum_{k=1}^{n-1}\bar{V}_{k},
$$
where and below by the sum over an empty set we mean zero.
It follows that, for $n\in[1,m]$, 
$$
L^{0}_{h}V_{n}\geq
2 \sum_{\lambda\in\Lambda^{n}}
( \bar\delta _{\lambda}u)
 \bar\delta _{\lambda}L^{0}_{h}u
+\delta\sum_{\lambda\in\Lambda^{n}}
\cQ( \bar\delta _{\lambda} u)
-K_{1}
\sum_{\lambda\in\Lambda^{n-1}}
\cQ( \bar\delta _{\lambda} u)
$$
\begin{equation}
                                             \label{4.8.1}
 -
(2c-2\delta c+\delta^{2} )
\bar{V}_{n}-N \sum_{k=1}^{n-1}\bar{V}_{k}.
\end{equation}
Next,
\begin{equation}                             \label{22.23.01.07}
2 \sum_{\lambda\in\Lambda^{n}}
( \bar\delta _{\lambda}u)
 \bar\delta _{\lambda}L^{0}_{h}u
=D_{t}V_{n}+R_{1}-R_{2},
\end{equation}
where
$$
R_{1}:=2 \sum_{\lambda\in\Lambda^{n}}
( \bar\delta _{\lambda}u)
 \bar\delta _{\lambda} (cu),
\quad 
R_{2}:=2 \sum_{\lambda\in\Lambda^{n}}
( \bar\delta _{\lambda}u)
 \bar\delta _{\lambda}f.
$$
Similarly to \eqref{4.8.3}  
$$
R_{1}=2cV_{n}+2\sum_{k=1}^{n}C^{k}_{n}R_{1k},
$$
where
$$
R_{1k}:=
\sum_{\lambda\in\Lambda^{n-k}}
\sum_{\mu\in\Lambda^{k}}( \bar\delta _{\mu} 
 \bar\delta 
_{\lambda}u)( \bar\delta _{\mu}c)
T_{\mu}( \bar\delta _{\lambda}u) .
$$
By our assumptions, 
\eqref{3.21.5}, and Cauchy's inequality 
$$
\sum_{\lambda\in\Lambda^{n-k}}
 |\bar\delta _{\mu} 
 \bar\delta 
_{\lambda}u|
T_{\mu}|\bar\delta _{\lambda}u|
\leq\big(\sum_{\lambda\in\Lambda^{n-k}}
 ( \bar\delta _{\mu} 
 \bar\delta 
_{\lambda}u) ^{2}\big)^{1/2}T_{\mu}
\big(\sum_{\lambda\in\Lambda^{n-k}}
 (   \bar\delta 
_{\lambda}u) ^{2}\big)^{1/2}
$$
$$
\leq\bar{V}_{n-k}^{1/2}\big(\sum_{\lambda\in\Lambda^{n-k}}
 ( \bar\delta _{\mu} 
 \bar\delta 
_{\lambda}u) ^{2}\big)^{1/2},
$$
$$
|R_{1k}|\leq\bar{V}_{n-k}^{1/2}V_{n}^{1/2}
\big(\sum_{\mu\in\Lambda^{k}}|\bar\delta_{\mu}c|^{2}\big)^{1/2}
\leq A^k
(\sup_{H_T}|D^kc|) V_{n}^{1/2}\bar{V}^{1/2}_{n-k}.
$$
 
Cauchy's inequality also allows us to estimate
$R_{2}$ and conclude from \eqref{4.8.1} that,
for $n\in[1,m]$,
$$
L^{0}_{h}V_{n}-2cV_{n}-D_{t}V_{n}\geq
 \delta \cQ_{n}
$$
\begin{equation}
                                             \label{4.8.2}
- K_{1}\cQ_{n-1}-
 (2c-2\delta c+2\delta^{2}) 
\bar{V}_{n}-N \sum_{k=0}^{n-1}\bar{V}_{k}-N
F_{n}^{2}, 
\end{equation}
 where  
$$
\cQ_{k}=\sum_{\lambda\in\Lambda^k}
\cQ(\bar\delta_{\lambda} u). 
$$
We now prove \eqref{4.8.03} by 
  showing that for each $n\in[1,m]$
\begin{equation}                            \label{8.19.01.08}
V_{k}\leq N (F_k^2+G_k^2) ,
\quad k=0,1,...,n.
\end{equation}
We prove this by induction on $n$. 
By Lemma \ref{lemma 3.6.1} we have 
$$
V_{0}\leq N (F_0^2+G_0^2). 
$$
Using this, from \eqref{4.8.2} we obtain 
$$
V_{1}\leq N (F_1^2+G_1^2)
$$ 
by Lemma \ref{lemma 3.6.1}, 
provided that 
$0\leq 2c-2\delta c+2\delta^{2}\leq 2c-\delta^{2}$
which is true indeed if
\begin{equation}
                                          \label{08.2.14.1}
3\delta\leq 2c_{0}.
\end{equation}
This may look like a nontrivial restriction on $\delta$.
However, obviously, 
if our assumptions are
satisfied with a $\delta\in(0,1)$, 
they are also satisfied with
any
$\delta'\in(0,\delta]$. Therefore, 
without losing generality
we suppose that \eqref{08.2.14.1} is valid. 
 Thus we have obtained 
\eqref{8.19.01.08} for $n=1$. Let $n\geq2$ and assume that 
\eqref{8.19.01.08} holds with $n-1$ in place of $n$. Then 
from \eqref{4.8.2}  
for $k=1,...,n$ we get  
\begin{equation}
                                          \label{08.2.14.2}
L^{0}_{h}V_{k}-2cV_{k}-D_{t}V_{k}\geq
 \delta \cQ_{k}- K_{1}\cQ_{k-1}-
 C_{\delta}I_{k=n}
\bar{V}_{n}-N(F_n^2+G_n^2) , 
\end{equation}
 with  
$C_{\delta}=2(c-\delta c+ \delta^{2})$. 
Actually, \eqref{08.2.14.2} is 
 true also  for $k=0$
if we set $\cQ_{-1}=0$, since
$$
L ^{0}_{h}(u^{2})-2cu^{2}-D_{t}(u^{2})
=2u(L ^{0}_{h}u-cu - 
D_{t}u)+ \cQ_{0}(u)
$$
$$
=-2uf+ \cQ_{0}(u) \geq-N(F_{1}^{2}+G_1^{2})
+ \cQ_{0}(u) .
$$
Next we set $\mu=K_{1}/\delta$, multiply \eqref{08.2.14.2}
by $\mu^{n-k}$ and sum up the resulting inequalities
 with respect to
$k=0,...,n$. Then, for
$$
W_{n}:=\sum_{k=0}^{n}\mu^{n-k}V_{k},
$$
we obtain
$$
L^{0}_{h}W_{n}-2cW_{n}-D_{t}W_{n}
\geq\delta\cQ_{n}-C_{\delta}\bar{V}_{n}
-N(F^{2}_{n}+G^{2}_{n})
$$
$$
\geq-C_{\delta}\bar{W}_{n}
-N(F^{2}_{n}+G^{2}_{n}).
$$
Recalling \eqref{08.2.14.1} and using Lemma \ref{lemma 3.6.1} 
shows that \eqref{8.19.01.08} holds.  
This justifies the induction and proves the theorem.

\mysection{Discussion of Assumptions     \protect\ref{assumption 11.22.11.06} 
and \protect\ref{assumption 07.10.16.1}}
                                               \label{section 4.9.1}

In \cite{GK} there are many sufficient conditions
for Assumption
\ref{assumption 11.22.11.06}  to be satisfied.
In this section   we suppose that only
Assumptions \ref{assumption 16.12.07.06} and
\ref{assumption 1.26.11.06} are satisfied and $m\geq2$. 
Assume also that for a number
$\bar\tau>  0 $ we 
have that, for any $\lambda\in\Lambda_{1}$,
\begin{equation}
                                            \label{08.2.16.1}
\text{either}
\quad\tau_{\lambda}\geq\bar{\tau}\quad
\text{or}\quad Dq_{\lambda}(t,x)=0\quad\text{for all}
\quad(t,x).
\end{equation}  
Recall that by $\cK$ we denote 
a generic operator from 
class $\frK$,   which may depend 
on $h$ and $t$, and may change each time it 
is mentioned even in one line.

\begin{remark}
                                   \label{remark 07.10.22.1}
Assume that
$m\geq2$, $\Lambda_{1}=-\Lambda_{1}$, $q_{\lambda}=
q_{-\lambda} $ 
and $p_{\lambda}\geq0$. Then, since $q_{\lambda}$ are
 twice continuously differentiable in $x$
and nonnegative by Assumption
 \ref{assumption 1.26.11.06}, we know that
$r_{\lambda}:=\sqrt{{q}_{\lambda}}$ is Lipschitz
continuous in $x$ with the Lipschitz constant
independent of $t$.

Conditions \eqref{3.24.1} and  
\eqref{3.24.01} involve a mixture
of finite differences and derivatives. 
Therefore, it is reasonable to 
try to find conditions
in terms only of finite differences which would imply
\eqref{3.24.1} and  \eqref{3.24.01}.
We claim that \eqref{3.24.1} and  \eqref{3.24.01}
are  satisfied with a $\tau_{0}>0$ and, perhaps,
different $\delta$, $\cK$, 
 $K_{1}$  if
for all smooth $\varphi$ on $H_{T}$
 and $n=1,...,m$ we have
$$
2m\sum_{\lambda\in\Lambda_{1}}
( \bar\delta _{ \lambda}\varphi)
L^{0}_{ \lambda}T_{ \lambda}\varphi\leq
(1- 
 \delta)\sum_{\lambda \in\Lambda_{1} }
\cQ( \bar\delta _{ \lambda}\varphi)
$$
\begin{equation}
                                           \label{4.9.3}
+K_{1}\cQ(\varphi)
+
(1-\delta)c\cK \big(\sum_{\lambda\in\Lambda_{1}}
| \bar\delta _{ \lambda}\varphi|^{2}\big) , 
\end{equation}
$$
2n\sum_{\lambda,\nu\in\Lambda_{1}}
( \bar\delta _{\lambda}
 \bar\delta _{\nu}\varphi)
L^{0}_{\lambda}T_{\lambda}
 \bar\delta _{\nu}\varphi
+
n(n-1)\sum_{\lambda\in\Lambda^{2}_{1}}
( \bar\delta _{\lambda}\varphi)
Q_{\lambda}T_{\lambda}\varphi
$$
$$
\leq (1- \delta)
\sum_{\nu\in\Lambda_{1}^{2}}
\cQ( \bar\delta _{\nu}\varphi)
+K_{1}\sum_{\nu\in\Lambda_{1}}
\cQ( \bar\delta _{\nu}\varphi)
$$
\begin{equation}
                                              \label{07.10.22.6}
+(1-\delta)c\cK\big(\sum_{\lambda\in\Lambda_{1}^{2} }
| \bar\delta _{\lambda}\varphi|^{2}\big)
+K_{1}\cK\big(\sum_{\lambda\in\Lambda_{1} }
| \bar\delta _{\lambda}\varphi|^{2}\big). 
\end{equation}
(Notice that the term $2(1-\delta)$
in \eqref{3.24.1} and  \eqref{3.24.01} is replaced now
with $1-\delta$.)

To prove this, first observe that 
as it follows from Remarks  5.1 
and 
5.2 of \cite{GK}, 
owing to \eqref{4.9.3}, 
the above mentioned properties of $r_{\lambda}$, and
\eqref{08.2.16.1}, there exist constants $\delta,
\bar{\tau}_{0}\in(
0,1]$, and $K_{1}$, perhaps different from the above ones,
  such that 
\begin{equation}
                                           \label{07.10.26.1}
m A_{1}(\varphi)
\leq(1- \delta)\sum_{\lambda\in\Lambda}
\cQ( \bar\delta _{\lambda}\varphi)
+K_{1}\cQ(\varphi)
+(1-\delta)c\cK\big(\sum_{\lambda\in\Lambda }
| \bar\delta _{\lambda}\varphi|^{2}\big) 
\end{equation}
on $H_{T}$
for all smooth functions $\varphi $  
provided that $\tau_{0}\in(0,\bar{\tau}_{0}]$. 
Thus, a condition even somewhat stronger than 
\eqref{3.24.1} is satisfied. 
Next, observe that 
the left-hand side of \eqref{3.24.01} equals 
$$
 B 
+A'+A''_{q}+A''_{p}+
B'+B'',
$$
 where 
 $B$ is the left-hand side of 
\eqref{07.10.22.6}, 
$$
A'=n\sum_{\nu\in\Lambda_{2}}
A_{1}( \bar\delta _{\nu}\varphi),
\quad
A''_{q}=2n\sum_{\lambda\in\Lambda_{1},\nu\in\Lambda_{2}}
 ( \bar\delta _{\nu}
 \bar\delta _{\lambda}
\varphi)Q_{\nu}
 \bar\delta _{\lambda}\varphi,
$$
$$
A''_{p}=2n\sum_{\lambda\in\Lambda_{1},\nu\in\Lambda_{2}}
 ( \bar\delta _{\nu}
 \bar\delta _{\lambda}
\varphi)P_{\nu}
 \bar\delta _{\lambda}\varphi,
$$
$$
B'=2n(n-1)\sum_{\lambda\in\Lambda_{1}\times\Lambda_{2}}
( \bar\delta _{\lambda}\varphi)
Q_{\lambda}T_{\lambda}\varphi,
$$
$$
B''=
n(n-1)\sum_{\lambda\in\Lambda_{2}^{2} }
( \bar\delta _{\lambda}\varphi)
Q_{\lambda}\varphi.
$$
Here  by \eqref{07.10.26.1}
$$
A'\leq(1- \delta)
\sum_{\nu\in\Lambda_{2},
\lambda\in\Lambda}
\cQ( \bar\delta _{\lambda}
 \bar\delta _{\nu}\varphi)
+K_{1}\sum_{\nu\in\Lambda_{ 2 }}
\cQ( \bar\delta _{\nu}\varphi)
$$
$$
+ (1-\delta)c\cK\big(\sum_{\lambda\in\Lambda^{2} }
| \bar\delta _{\lambda}\varphi|^{2}\big) .
$$
Then,
$$
A''_{q}=2n
 \tau_0^2 
\sum_{\lambda,\mu\in\Lambda_{1}}
\sum_{j=1}^{d}(D_{j}
 \bar\delta _{\lambda}\varphi) 
(D_{j}r_{\mu})[r_{\mu}
\Delta_{\mu}
 \bar\delta _{\lambda}\varphi]
$$
$$
\leq (1/16)\delta c
 \tau_0^2 \sum_{\lambda,\mu\in\Lambda_{1}}
\sum_{j=1}^{d}(D_{j}
 \bar\delta _{\lambda}\varphi)^{2}
+N
 \tau_0^2 
\sum_{\mu\in\Lambda_{1},\lambda\in\Lambda_{1}^{2}}q_{\mu}
( \bar\delta _{\mu}
 \bar\delta _{\lambda}\varphi)^{2}
$$
$$
=(1/16)\delta c\sum_{\lambda\in\Lambda_{1},
\nu\in\Lambda_{2}}
( \bar\delta _{\nu}
 \bar\delta _{\lambda}\varphi)^{2}
+N
 \tau_0^2 
\sum_{\mu\in\Lambda_{1},\lambda\in\Lambda_{1}^{2}}q_{\mu}
( \bar\delta _{\mu}
 \bar\delta _{\lambda}\varphi)^{2}
$$
$$
\leq(1/16)\delta c\sum_{\lambda\in\Lambda^{2}}
( \bar\delta _{\lambda}\varphi)^{2}
+N
 \tau_0^2 
 \sum_{\lambda\in\Lambda^{2}}
\cQ( \bar\delta _{\lambda}\varphi),
$$
where and below by $N$ we denote various generic constants
independent of $\varphi$, $(t,x)$, and 
$\tau_{0}$. 
Next, quite similarly
$$
A''_{p}=2n
 \tau^2_0 
\sum_{\lambda,\mu\in\Lambda_{1}}
\sum_{j=1}^{d}(D_{j}
 \bar\delta _{\lambda}\varphi)
(D_{j}p_{\mu}) 
 \bar\delta _{\mu}
 \bar\delta _{\lambda}\varphi
$$
$$
\leq(1/16)\delta c
 \tau_0^2  \sum_{\lambda\in\Lambda_{1}}
\sum_{j=1}^{d}
(D_{j} \bar\delta _{\lambda}\varphi)^{2}
+N
 \tau_0^2 
\sum_{\lambda\in\Lambda_{1}^{2}}
( \bar\delta _{\lambda}\varphi)^{2}
$$
$$
\leq(1/16)\delta c\sum_{\lambda\in\Lambda^{2}}
( \bar\delta _{\lambda}\varphi)^{2}
+N
 \tau_0^2 \sum_{\lambda\in\Lambda^{2}}
( \bar\delta _{\lambda}\varphi)^{2}.
$$
Now we estimate $B'$ and $B''$. We have
$$
B'=
n(n-1)
 \tau^2_0 \sum_{\lambda,\mu\in\Lambda_{1} }
\sum_{j=1}^{d}
(D_{j} \bar\delta _{\lambda}\varphi)
(D_{j} \bar\delta _{\lambda}q_{\mu})
\Delta_{\mu}
T_{\lambda}\varphi.
$$
Here
$$
\Delta_{\mu}T_{\lambda}=\Delta_{\mu}
+(\delta_{\mu}+\delta_{-\mu})\delta_{\lambda}, 
$$
and it is seen that
$$
B'\leq(1/16)\delta c
 \tau_0^2 
\sum_{\lambda\in\Lambda_{1}}\sum_{j=1}
^{d}
(D_{j}
 \bar\delta _{\lambda}\varphi)^{2}
+N
 \tau_0^2 
\sum_{\lambda\in\Lambda^{2}}
( \bar\delta _{\lambda}
\varphi)^{2}
$$
$$
\leq(1/16)\delta c\sum_{\lambda\in\Lambda^{2}}
( \bar\delta _{\lambda}\varphi)^{2}
+N
 \tau_0^2 
\sum_{\lambda\in\Lambda^{2}}
( \bar\delta _{\lambda}
\varphi)^{2}.
$$
Similarly,
$$
B''=
(1/2)
n(n-1)
 \tau_0^4 
\sum_{j,k=1 }^{d}\sum_{\mu\in\Lambda_{1}}
(D_{jk}\varphi)
(D_{jk}q_{ \mu})\Delta_{\mu}\varphi
$$
$$
\leq(1/16)\delta c\sum_{\lambda\in \Lambda^{2}}
( \bar\delta _{\lambda} \varphi)^{2}
+N
 \tau_0^4 
\sum_{\lambda\in\Lambda^{2}}
( \bar\delta _{\lambda}
\varphi)^{2}.
$$

By combining the above estimates we see that
the left-hand side of \eqref{3.24.01} is majorated by
$$
 (1- \delta+N
 \tau_0^2 )\sum_{ \nu\in\Lambda^{2}}
\cQ(
 \bar\delta _{\nu}\varphi)+ 
K_{1}\sum_{ \nu\in\Lambda }
\cQ(
 \bar\delta _{\nu}\varphi)+K_{1}
\cK\big(\sum_{\lambda\in\Lambda }
| \bar\delta _{\lambda}\varphi|^{2}\big) 
$$
$$
+ c\big[2(1-\delta)\cK\big(\sum_{\lambda\in\Lambda^{2} }
| \bar\delta _{\lambda}\varphi|^{2}\big) 
+((1/4)\delta+N
 \tau_0^2 )\sum_{\lambda\in\Lambda^{2} }
| \bar\delta _{\lambda}\varphi|^{2}\big].
$$
 It follows easily that by choosing 
 $\tau_0$ small  enough
we will satisfy \eqref{3.24.01}
as well as \eqref{3.24.1} 
with $\delta/2$ in place of $\delta$
and appropriate $K_{1}$.

\end{remark}
\begin{remark}
                                   \label{remark 07.10.17.1}

In \cite{GK} we have seen that even 
Assumption
\ref{assumption 11.22.11.06} imposes certain nontrivial 
{\em structural\/} conditions on $q_{\lambda}$ which cannot
be guaranteed by  
 the size of $c_{0}$ if  $q_{\lambda}$ is only once 
continuously differentiable.
In contrast, given that Assumptions  
\ref{assumption 16.12.07.06},
\ref{assumption 18.12.07.06},  
\ref{assumption 11.22.11.06}
are satisfied and $m\geq2$, 
we claim that Assumption \ref{assumption 07.10.16.1}
is also satisfied if $c_{0}$ is large enough.
 
 To prove our claim we notice that
by \eqref{3.22.3} 
$$
\sum_{\lambda\in\Lambda^{2}} |
  Q_{\lambda}T_{\lambda}\varphi|^{2}
\leq  
 N
\cK 
(\sum_{ \mu\in\Lambda^{ 2}} 
| \bar\delta _{\mu} \varphi |^{2}) ,
$$
so that
$$
n(n-1)\sum_{\lambda\in\Lambda^{2}} (
 \bar\delta _{\lambda}
\varphi)
  Q_{\lambda}T_{\lambda}\varphi
$$
$$
\leq m(m-1)\big[\sum_{\lambda\in\Lambda^{2}}
| \bar\delta _{\lambda}\varphi|^{2}+
 N\cK (\sum_{ \mu\in\Lambda^{ 2}} 
| \bar\delta _{\mu} \varphi |^{2})\big]=:
N'\cK (\sum_{ \mu\in\Lambda^{ 2}} 
| \bar\delta _{\mu} \varphi |^{2}).
$$
Now assume that $c$ is so large that
$$
N' \leq\delta c.
$$
Then it follows from \eqref{3.24.1} that
the left-hand side of \eqref{3.24.01} is majorated by
$$
(1-\delta)\sum_{\lambda \in\Lambda_{1},\nu\in\Lambda^{2}}
\chi_{\lambda}|\delta_{\lambda}
 \bar\delta _{\nu}\varphi|^{2}
+K_{1}\sum_{\lambda\in\Lambda_{1},\nu\in\Lambda}\chi_{\lambda}
|\delta_{\lambda}
 \bar\delta _{\nu}\varphi|^{2}+I,
$$
where
$$
I=2(1-\delta)c\cK\big(\sum_{\lambda\in\Lambda^{2} }
| \bar\delta _{\lambda}\varphi|^{2}\big) 
+
\delta c
 \cK \big(\sum_{ \mu\in\Lambda^{ 2}} 
| \bar\delta _{\mu} \varphi |^{2}\big)
$$
$$
=:2(1-\delta/2)c \cK 
\big(\sum_{\lambda\in\Lambda^{2} }
| \delta _{\lambda}\varphi|^{2}\big).
$$
  We thus obtain \eqref{3.24.01}
with $\delta/2$ in place of $\delta$.
\end{remark}

\begin{remark}
It is interesting to have sufficiently simple conditions
on the coefficients of 
 differential operators  $\cL$ which guarantee that there 
exist finite-difference schemes for which our
assumptions hold. Here we will only give 
a one dimensional
example. This example is based on the results of Remark
\ref{remark 07.9.18.5} below, which can also be used
to analyze many multi-dimensional situations as well
in the spirit of the comments in \cite{GK}.

Take $d=1$ and
$$
\cL\varphi(x)=(1/2)a(x)\varphi''(x)+b(x)
\varphi'(x)-c(x)\varphi(x).
$$
We assume that $a\geq0$ and $r:=\sqrt{a},b$, and $c$
are $m$-times 
continuously differentiable   with bounded 
derivatives. 
We take $\Lambda_{1}=\{\pm1\}$ and define
$$
q_{\mu}=a, \quad p_{\mu}=(1/2)\mu b+\theta,
$$
where $\theta$  is a   constant such that $p_{\mu}\geq1$.
By using an argument  in Remark 6.7 of \cite{GK}
and using our Remark \ref{remark 07.9.18.5}, 
one can easily derive that, for a sufficiently
small $\tau_{0}$ and $\tau_{\mu}\equiv1$, 
Assumptions \ref{assumption 11.22.11.06}
and \ref{assumption 07.10.16.1}
 are satisfied for all sufficiently small $h$
(with perhaps different $\delta$ and $K_{1}$)
if, for $n\leq m$,
$$
75n^{2}(r')^{2}+2nb'\leq(1-\delta)c+K_{1}a.
$$
Again as in \cite{GK} we see that at points
where $a$ is close to zero either $c$ should be large
or $b'$ be sufficiently negative. 
\end{remark}

\begin{remark}
                                      \label{remark 07.9.18.5}
Condition \eqref{4.9.3} and its implications are  discussed
in many details
 in \cite{GK} (with $2c$ in place of $c$).
Here we give sufficient conditions for \eqref{4.9.3} and
\eqref{07.10.22.6}
to be satisfied
 without involving test functions $\varphi$. 
 For simplicity,  we only do it in case 
$$
\tau_{\lambda}=1\quad\text{for all $\lambda\in\Lambda_1$}.
$$
It is obvious that if we define $\xi_{\lambda\mu}=\delta_{\lambda}
\delta_{\mu}\varphi$, then condition
\eqref{07.10.22.6} can be rewritten in terms of $\xi_{\lambda\mu}$.
What is nontrivial is that one can give sufficient
conditions for \eqref{07.10.22.6} to hold in terms
of $\xi_{\lambda}$ and not the
two-parameter object
$\xi_{\lambda\mu}$. In addition,
we will see that  these  sufficient conditions
are obtained just by {\em slightly strengthening\/}
 the corresponding conditions from
\cite{GK} guaranteeing the first-order derivatives estimates.
As in \cite{GK} one could extract further implications
and simplifications of the new conditions of the type that
on the set where $c$ is small we need $ \chi _{\lambda}$ to be
uniformly bounded away from zero or $p_{\lambda}$ be
sufficiently strongly monotone (see \cite{GK} for 
more details).

 As in Remark \ref{remark 07.10.22.1}
we assume that $\Lambda_{1}=-\Lambda_{1}$ and $q_{\lambda}=
q_{-\lambda}$ ($\geq0$), $p_{\lambda}\geq0$,  and, 
additionally,
$r_{\lambda}:=\sqrt{{q}_{\lambda}}$ is {\em twice\/} continuously
differentiable in $x$ and is bounded on $H_{T}$
along with first and second-order derivatives in $x$. 
Also we fix a constant  $\delta\in(0,1/ 4 ]$ and
assume that on $H_{T}$
there are functions $r_{\lambda\mu}=r_{h\lambda\mu}$, 
$p_{\lambda\mu}=p_{h\lambda\mu}\geq0$,
$\lambda,\mu\in\Lambda_{1}$,
such that 
\begin{equation}
                                                 \label{07.9.17.2}
m(m-1)h^{2}(\delta_{\lambda}r_{\mu})^{2}
\leq  \delta(
\chi_{\lambda}+\chi_{\mu})
+h^{2}r^{2}_{\lambda\mu},\quad\sum_{\mu\in\Lambda_{1}}
\sup_{\lambda\in\Lambda_{1}}r^{2}_{\lambda\mu}\leq\delta c,
\end{equation}
\begin{equation}
                                                 \label{07.9.17.3}
h^{2}|\delta_{\lambda}p_{\mu}|\leq\delta^{2}
(
\chi_{\lambda}+\chi_{\mu})
+\delta h^{2}p_{\lambda\mu},
\quad\sum_{\mu\in\Lambda_{1}}
\sup_{\lambda\in\Lambda_{1}}p_{\lambda\mu}\leq\delta c.
\end{equation}
By virtue of Remark  6.1  of \cite{GK} 
one can always find approximations $L'_{h}$
 of the zero operator such that $L_{h}+L'_{h}$
will still be approximating $\cL$ and for the coefficients
$p'_{\lambda}$ of $L_{h}+L'_{h}$ we will 
have $p'_{\lambda}\geq1$. Obviously, 
for $L_{h}+L'_{h}$ conditions
\eqref{07.9.17.2} and \eqref{07.9.17.3} are satisfied
with $r_{\lambda\mu}=p_{\lambda\mu}=0$ 
for sufficiently small $h$.

For a function $\xi_{\lambda}$ given on $\Lambda_{1}$ let us write
$$
|\xi|^{2}=\sum_{\lambda\in\Lambda_{1}}|\xi_{\lambda}|^{2}
$$
and let us drop the summation sign over repeated indices
in $\Lambda_{1}$.
Then we claim that conditions \eqref{4.9.3} and
\eqref{07.10.22.6}
are satisfied with appropriate $K_{1}$, 
 $\cK$,
and $\delta$
if on $H_{T}$ for  all 
functions $\xi_{\lambda}$  and $n=1,...,m$
we have
$$
28 n^{2}
(1- 4 \delta)^{-1}J_{1}
+(9/2)n^{2}(1-4\delta)^{-1}J_{2}+(1/2)n^{2}J_{3}
$$
$$
+2\delta n^{2}
  \sum_{\lambda,\mu \in\Lambda_{1}} \xi_{\lambda} ^{2}
|\delta_{\lambda}p_{\mu}|  +2n\xi_{\lambda}\xi_{\mu}
 (\delta_{\lambda}
p_{\mu}
+
(\delta_{\lambda}r_{\mu})^2) 
$$
\begin{equation}
                                             \label{07.9.16.1}
\leq(1 -4\delta)c|\xi|^{2}
+K_{1} 
 \xi_{\lambda} ^{2}
 \chi _{\lambda}
+ \delta h^{-2} 
 \chi _{\lambda}
|\xi_{\lambda}+\xi_{-\lambda}|^{2},
\end{equation}
where
$$
J_{1}=\sum_{\mu,\lambda\in\Lambda_{1}}
\xi_{\lambda}^{2}(\delta_{\lambda}r_{\mu})^{2},\quad
J_{2}=\sum_{\mu\in\Lambda_{1}}\big(
\sum_{\lambda\in\Lambda_{1}}\xi_{\lambda}\delta_{\lambda}
r_{\mu}\big)^{2},\quad
J_{3}=\sum_{\lambda,\mu\in\Lambda_{1}}
(\delta_{\lambda}r_{\mu})^{2}
\xi_{\mu}^{2}.
$$
 
To prove this claim, introduce  
$$
J_{1}(\varphi)=\sum_{\mu,\lambda\in\Lambda_{1}}
(\delta_{\lambda}\varphi)^{2}(\delta_{\lambda}r_{\mu})^{2},
$$
$$
J_{2}(\varphi)=\sum_{\mu\in\Lambda_{1}}\big(
\sum_{\lambda\in\Lambda_{1}}(\delta_{\lambda}\varphi)
\delta_{\lambda}
r_{\mu}\big)^{2},
\quad
J_{3}(\varphi)=\sum_{\lambda,\mu\in\Lambda_{1}}(
\delta_{\mu}\varphi)^{2}(\delta_{\lambda}r_{\mu})^{2}
$$
and
first recall that by Remarks 
 6.2 and 6.3  of \cite{GK} 
after replacing there $c$ with $c/2$ we obtain
$$
2n
 \sum_{\lambda\in\Lambda_1} 
(\delta_{\lambda}\varphi)L^{0}_{\lambda}T_{\lambda}\varphi
+
18
n^{2}J_{1}(\varphi)+(5/2)n^{2}J_{2}(\varphi)
+(1/2)n^{2}J_{3}(\varphi)
$$
$$
\leq (1-\delta)\sum_{\lambda\in\Lambda_{1}}\cQ(\delta_{\lambda}
\varphi)+K_{1}
\cQ(\varphi)+(1-\delta )c\cK\big(\sum_{\lambda\in
\Lambda_{1}}|\delta_{\lambda}\varphi|^{2}\big).
$$
 In particular condition
\eqref{4.9.3} is satisfied. Furthermore,
by substituting $\delta_{\nu}\varphi$ in place of $\varphi$
and summing up over $\nu\in\Lambda_{1}$, we get
$$
2n\sum_{\lambda,\nu\in\Lambda_{1}}(\delta_{\nu}
\delta_{\lambda}\varphi)L^{0}_{\lambda}T_{\lambda}
\delta_{\nu}\varphi
+
18
n^{2}\sum_{\nu\in\Lambda_{1}}
J_{1}(\delta_{\nu}\varphi)
$$
$$
+(5/2)n^{2}
\sum_{\nu\in\Lambda_{1}}J_{2}(\delta_{\nu}\varphi)
+(1/2)
n^{2}\sum_{\nu\in\Lambda_{1}}J_{3}(\delta_{\nu}\varphi)
$$
$$
\leq (1-\delta)\sum_{\lambda\in\Lambda_{1}^{2}}\cQ(\delta_{\lambda}
\varphi)+K_{1}
\sum_{\nu\in\Lambda_{1}}\cQ(\delta_{\nu}\varphi)
+(1-\delta )c\cK\big(\sum_{\lambda\in
\Lambda_{1}^{2}}|\delta_{\lambda}\varphi|^{2}\big).
$$

It follows that to prove our claim, it suffices to prove
that
$$
n(n-1)\sum_{\lambda\in\Lambda_{1}^{2}}
(\delta_{\lambda}\varphi)Q_{\lambda}T_{\lambda}
\varphi\leq
18
n^{2}\sum_{\nu\in\Lambda_{1}}
J_{1}(\delta_{\nu}\varphi)+(5/2)n^{2}
\sum_{\nu\in\Lambda_{1}}J_{2}(\delta_{\nu}\varphi)
$$
$$
+(1/2)
n^{2}\sum_{\nu\in\Lambda_{1}}J_{3}(\delta_{\nu}\varphi)+
(2/3)\delta\sum_{\nu\in\Lambda_{1}^{2}}\cQ(\delta_{\nu}
\varphi)
$$
\begin{equation}
                                             \label{10.10.02.08}
+(1/3)\delta c\cK\big(\sum_{\nu\in\Lambda_{1}^{2}}
|\delta_{\nu}\varphi|^{2}\big)
+N\sum_{\nu\in\Lambda_{1}}\cQ(\delta_{\nu}\varphi)
+N\cK\big(\sum_{\nu\in\Lambda_{1}}
|\delta_{\nu}\varphi|^{2}\big),
\end{equation}
where and below by $N$ we denote generic constants
independent of $\varphi$ and $(t,x)$ (and 
various $\varepsilon $'s
once they appear).
Observe that for 
$\lambda
 =(\lambda^{1},\lambda^{2}) 
\in\Lambda_{1}^{2}$
and $\mu\in\Lambda_{1}$,
 $$
T_{\lambda}
=1+h^2\delta_{\lambda^1}\delta_{\lambda^2}
+h(\delta_{\lambda^1}+\delta_{\lambda^2}), 
\quad
\Delta_{\mu}=h^{-1}(\delta_{\mu}+\delta_{-\mu})
$$
and hence 
$$
\Delta_{\mu}T_{\lambda}=\Delta_{\mu}
+h\delta_{\lambda^{1}}\delta_{\lambda^{2}}(\delta_{\mu}
+\delta_{-\mu})+(
\delta_{\lambda^{1}}
+\delta_{\lambda^{2}})(\delta_{\mu}
+\delta_{-\mu}),
$$
implying that 
$$
\sum_{\lambda\in\Lambda_{1}^{2}}(\delta_{\lambda}\varphi)
Q_{\lambda}T_{\lambda}\varphi= S_{1}+S_{2},
$$
where
$$
S_{1}=
(1/2)
\sum_{\lambda\in\Lambda_{1}^{2},\mu\in\Lambda_1} 
(\delta_{\lambda}\varphi)
(\delta_{\lambda}q_{\mu}) (4
\delta_{\lambda^{1}}-\delta_{-\mu})\delta_{\mu}\varphi,
$$
$$
S_{2}= 
h\sum_{\lambda\in\Lambda_{1}^{2},\mu\in\Lambda_1} 
(\delta_{\lambda}\varphi)
(\delta_{\lambda}q_{\mu})
\delta_{\lambda} \delta_{\mu}
\varphi.
$$
Next, as it is easy to see for $\lambda\in\Lambda_{1}^{2}$
$$
\delta_{\lambda}q_{\mu}=2(\delta_{\lambda^{1}}r_{\mu})
\delta_{\lambda^{2}}r_{\mu}
+
2r_{\mu}\delta_{\lambda}r_{\mu}
$$
$$
+2h(\delta_{\lambda^{1}}r_{\mu}+\delta_{\lambda^{2}}r_{\mu})
\delta_{\lambda}r_{\mu}+h^{2}(\delta_{\lambda}r_{\mu})^{2}.
$$

{\em Estimating $S_{2}$\/}.
First we estimate the term $S_{2}$, which contains the third-order
differences of $\varphi$. For the main term in $S_{2}$ we have
$$
B_{1}:=2
h\sum_{\lambda\in\Lambda_{1}^{2}}(\delta_{\lambda}\varphi)
(\delta_{\lambda^{1}}r_{\mu})(
\delta_{\lambda^{2}}r_{\mu})
\delta_{\lambda} \delta_{\mu}
\varphi
$$
$$
\leq
 16 
\sum_{\lambda,\mu,\nu\in\Lambda_{1}}
(\delta_{\lambda}\delta_{\nu}\varphi)^{2}(\delta_{\lambda}
r_{\mu})^{2}
+(1/16)h^{2}\sum_{\lambda,\mu,\nu\in\Lambda_{1}}
(\delta_{\lambda}r_{\mu})^{2}(\delta_{\lambda}
\delta_{\mu}\delta_{\nu}\varphi)^{2}
$$
$$
=:16 
\sum_{\nu\in\Lambda_{1}}J_{1}(\delta_{\nu}
\varphi)+(1/16)E,
$$
where by assumption \eqref{07.9.17.2} 
and Lemma  6.1 of \cite{GK}
$$
n(n-1)
E\leq2\delta\sum_{\nu\in\Lambda_{1}^{2}}\cQ(\delta_{\nu}
\varphi)+4\delta c \cK\big(\sum_{\lambda\in\Lambda_{1}^{2}}
(\delta_{\lambda}\varphi)^{2}\big).
$$
Hence,
$$
n(n-1)
B_{1}
\leq 
 16 
n^{2}
\sum_{\nu\in\Lambda_{1}}J_{1}(\delta_{\nu}
\varphi)
$$
\begin{equation}
                                            \label{10.22.1}
+(1/8)\delta
\sum_{\nu\in\Lambda_{1}^{2}}\cQ(\delta_{\nu}
\varphi)+
(1/4)
\delta c \cK\big(\sum_{\lambda\in\Lambda_{1}^{2}}
(\delta_{\lambda}\varphi)^{2}\big).
\end{equation}
Next, obviously, for any $\varepsilon>0$, (here
we use that $p_{\lambda}\geq0$)
$$
B_{2}:=2
h\sum_{\lambda\in\Lambda_{1}^{2} ,\mu\in\Lambda_1} 
(\delta_{\lambda}\varphi)
r_{\mu} (
\delta_{\lambda}r_{\mu})
\delta_{\lambda} \delta_{\mu}
\varphi
$$
$$
\leq\varepsilon^{-1}
h^{2}
\sum_{\lambda\in\Lambda_{1}^{2},\mu\in\Lambda_{1}}
(\delta_{\lambda}\varphi)^{2}
(
\delta_{\lambda}r_{\mu})^{2}
+ \varepsilon
\sum_{\lambda\in\Lambda_{1}^{2},\mu\in\Lambda_{1}}q_{\mu}
(\delta_{\lambda} \delta_{\mu}
\varphi)^{2}
$$
$$
\leq N\varepsilon^{-1}
\cK\big(\sum_{\lambda\in\Lambda_{1}}
(\delta_{\lambda}\varphi)^{2}\big)
+ \varepsilon\sum_{\nu\in\Lambda_{1}^{2}}
\cQ(\delta_{\nu}\varphi).
$$
It follows that (with $\varepsilon>0$ different from
the one from above but still arbitrary)
\begin{equation}
                                            \label{10.22.2}
n(n-1)
B_{2} 
\leq N \varepsilon^{-1}
\cK\big(\sum_{\lambda\in\Lambda_{1}}
(\delta_{\lambda}\varphi)^{2}\big)
+\varepsilon\sum_{\nu\in\Lambda_{1}^{2}}
\cQ(\delta_{\nu}\varphi).
\end{equation}

 Also
$$
n(n-1)h^{2}\sum_{\lambda\in\Lambda_{1}^{2},
\mu\in\Lambda_{1}}
|(\delta_{\lambda}\varphi)
\delta_{\lambda} \delta_{\mu}
\varphi|
$$
$$
\leq\varepsilon\delta c
\sum_{\lambda\in\Lambda_{1}^{2}}(\delta_{\lambda}\varphi)^{2}+
N\varepsilon^{-1}\cK\big(\sum_{\lambda\in\Lambda_{1}}
(\delta_{\lambda}\varphi)^{2}\big).
$$

Upon combining this with \eqref{10.22.1} and \eqref{10.22.2}
we obtain
$$
n(n-1)S_{2}\leq(\varepsilon+\delta/8)
\sum_{\nu\in\Lambda_{1}^{2}}
\cQ(\delta_{ \nu}\varphi)+
 16 n^{2}
\sum_{ \nu \in\Lambda_{1}}
J_{1}(\delta_{\nu}\varphi)
$$
\begin{equation}
                                            \label{10.22.3}
+(\varepsilon+ 1/4)\delta c\cK\big(
\sum_{\lambda\in\Lambda_{1}^{2}}(\delta_{\lambda}\varphi)^{2}
\big)
+N\varepsilon^{-1}\cK\big(\sum_{\lambda\in\Lambda_{1}}
(\delta_{ \lambda}\varphi)^{2}\big).
\end{equation}

{\em Estimating $S_{1}$\/}. We again start with the main term in
$S_{1}$, which we split into two parts writing  
$$
(4
\delta_{\lambda^{1}}-\delta_{-\mu})\delta_{\mu}\varphi 
=4\delta_{\lambda^{1}}\delta_{\mu}\varphi
+ \Delta_{\mu}\varphi.
$$ 
We have 
$$
4\sum_{\lambda\in\Lambda_{1}^{2},\mu\in\Lambda_{1}}
(\delta_{\lambda}\varphi)
(\delta_{\lambda^{2}}
r_{\mu})(\delta_{\lambda^{1}}r_{\mu})\delta_{\lambda^{1}}
\delta_{\mu}\varphi
$$
$$
=4\sum_{\nu,\mu\in\Lambda_{1}}
\big[\sum_{\lambda\in\Lambda_{1}}
(\delta_{\lambda}(\delta_{\nu}\varphi)
) \delta_{\lambda}r_{\mu} \big]
(\delta_{\nu}r_{\mu})\delta_{\nu}\delta_{\mu}\varphi
$$
$$
\leq2\sum_{\nu,\mu\in\Lambda_{1}}\big[
\sum_{\lambda\in\Lambda_{1}}
(\delta_{\lambda}\delta_{\nu}\varphi)
 \delta_{\lambda}r_{\mu}\big]^{2}
+2\sum_{\nu,\mu\in\Lambda_{1}}(\delta_{\nu}r_{\mu})^{2}
(\delta_{\nu}\delta_{\mu}\varphi)^{2}
$$
$$
\leq2\sum_{\nu\in\Lambda_{1}}J_{2}(\delta_{\nu}\varphi)
+2\sum_{\nu,\mu,\lambda\in\Lambda_{1}}(\delta_{\nu}r_{\mu})^{2}
(\delta_{\nu}\delta_{\lambda}\varphi)^{2}
$$
$$
=2\sum_{\nu\in\Lambda_{1}}J_{2}(\delta_{\nu}\varphi)
+2\sum_{\lambda\in\Lambda_{1}}J_{1}(\delta_{\lambda}\varphi).
$$
Furthermore,
$$
\sum_{\lambda\in\Lambda_{1}^{2}
,\mu\in\Lambda_{1}}(\delta_{\lambda}\varphi)
(
\delta_{\lambda^{1}}r_{\mu}) (\delta_{\lambda^{2}}r_{\mu})
 \Delta_{\mu}\varphi
$$
$$
=
\sum_{\nu,\mu\in\Lambda_{1}}\big[\sum_{\lambda\in\Lambda_{1}}
(\delta_{\lambda}\delta_{\nu}\varphi) 
\delta_{\lambda}r_{\mu} \big](\delta_{\nu}r_{\mu})
\Delta_{\mu}\varphi
$$
$$
\leq
(1/2)
\sum_{\nu,\mu\in\Lambda_{1}}\big[\sum_{\lambda\in\Lambda_{1}}
(\delta_{\lambda}\delta_{\nu}\varphi)
 \delta_{\lambda}r_{\mu} \big]^{2}
+
(1/2)
\sum_{\nu,\mu\in\Lambda_{1}}
(\delta_{\nu}r_{\mu})^{2}
(\Delta_{\mu}\varphi)^{2}
$$
$$
\leq
(1/2)
\sum_{\nu\in\Lambda_{1}}J_{2}(\delta_{\nu}\varphi)
+
(1/2)
\sum_{\nu\in\Lambda_{1}}J_{3}(\delta_{\nu}\varphi).
$$

 Next, obviously
$$
4
\sum_{\lambda\in\Lambda_{1}^{2}
 \mu\in\Lambda_1} 
[(\delta_{\lambda}\varphi)
 (\delta_{\lambda}r_{\mu}) ]r_{\mu}
\delta_{\lambda^{1}}\delta_{\mu}\varphi
\leq\varepsilon c\delta
\sum_{\lambda\in\Lambda_{1}^{2}}(\delta_{\lambda}\varphi)^{2}
+\varepsilon^{-1}N
 \sum_{\lambda\in\Lambda_{1}}
\cQ(\delta_{\lambda}\varphi) ,
$$
$$
\sum_{\lambda\in\Lambda_{1}^{2}}[(\delta_{\lambda}\varphi)
 (\delta_{\lambda}r_{\mu}) ]r_{\mu}
\Delta_{\mu}\varphi
\leq\varepsilon c\delta
\sum_{\lambda\in\Lambda_{1}^{2}}(\delta_{\lambda}\varphi)^{2}
+\varepsilon^{-1}N 
 \sum_{\lambda\in\Lambda_{1}}
\cQ(\delta_{\lambda}\varphi) .
$$
Finally,
$$
h\sum_{\lambda\in\Lambda_{1}^{2},\mu\in\Lambda_{1}}
|(\delta_{\lambda}\varphi)
\delta_{\lambda^{1}}\delta_{\mu}\varphi|
\leq
\varepsilon c\delta
\sum_{\lambda\in\Lambda_{1}^{2}}(\delta_{\lambda}\varphi)^{2}
+\varepsilon^{-1}N \cK
\big(\sum_{\lambda\in\Lambda_{1}}(\delta_{\lambda}\varphi)^{2}
\big).
$$
Upon combining the above estimates we obtain
$$
n(n-1)S_{1}\leq n^{2}
\sum_{\nu\in\Lambda_{1}}
[2J_{1}(\delta_{\nu}\varphi)
+(5/2)J_{2}(\delta_{\nu}\varphi)
+(1/2)J_{3}(\delta_{\nu}\varphi)]
$$
$$
+\varepsilon c\delta
\sum_{\lambda\in\Lambda_{1}^{2}}(\delta_{\lambda}\varphi)^{2}
+\varepsilon^{-1}N\bigg(\sum_{\lambda\in\Lambda_{1}}
\cQ(\delta_{\lambda}\varphi)
+\cK
\big(\sum_{\lambda\in\Lambda_{1}}(\delta_{\lambda}\varphi)^{2}
\big)\bigg). 
$$
This along with \eqref{10.22.3}
leads to \eqref{10.10.02.08} after appropriately choosing
  $\varepsilon$ and proves our claim. 
\end{remark}

\end{document}